\documentclass[11pt]{article}

\usepackage{float}
\usepackage{amsthm,,amssymb, color, epsfig, graphics}
\usepackage{xcolor}
\usepackage{graphicx}
\usepackage[intlimits]{amsmath}
\usepackage{latexsym}

\newtheorem{theorem}{Theorem}

\newtheorem{proposition}{Proposition}

\theoremstyle{definition}
\newtheorem{definition}{Definition}
\newtheorem{example}{Example} 

\newtheorem{remark}{Remark}

\newcommand{\beq}{\begin{equation}}  
\newcommand{\eeq}{\end{equation}}  
\newcommand{\bea}{\begin{eqnarray}}  
\newcommand{\eea}{\end{eqnarray}}  
\newcommand\la{{\lambda}}

\newcommand\al{{\alpha}}   
\newcommand\bet{{\beta}} 
\newcommand\ze{{\zeta}} 
\newcommand\gam{{\gamma}}     
\newcommand\om{{\omega}}

\newcommand\rd{{\mathrm{d}}}  

\newcommand\ri{{\mathrm{i}}}

\newcommand{\Q}{{\mathbb Q}}
\newcommand{\Z}{{\mathbb Z}}
\newcommand{\C}{{\mathbb C}}
\newcommand{\R}{{\mathbb R}}

\newcommand{\T}{{\mathbb T}}
\newcommand{\Ps}{{\mathbb P}}

\newcommand\bx{{\bf x}}
\newcommand\bd{{\bf d}}

\newcommand\SH{{\mathcal{S}}}

\begin{document}

%

%

\title{Cluster algebras and discrete integrability}

\author{Andrew N.W. Hone\footnote{Currently on leave at UNSW,   Sydney, Australia.}, Philipp Lampe and Theodoros E. Kouloukas \\
School of Mathematics, Statistics \& Actuarial Science\\ 
University of Kent, Canterbury CT2 7FS, UK.}

\maketitle

\begin{abstract}
  \noindent
Cluster algebras are a class of commutative algebras whose generators are defined by a recursive process called mutation.
 We give a brief introduction to cluster algebras, and explain how discrete integrable systems can appear in the context 
of cluster mutation. In particular, we give examples of birational maps that are integrable in the Liouville sense and arise 
from cluster algebras with periodicity, as well as examples of discrete Painlev\'e equations that are derived from Y-systems. 
\end{abstract}

\section{Introduction}

Cluster algebras are a special class of commutative algebras that were introduced by 
Fomin and Zelevinsky almost twenty years ago \cite{fz1}, and rapidly became the hottest topic in modern algebra. 
Rather than being defined a priori by a given set of 
generators and relations, the generators of a cluster algebra are produced recursively by 
iteration of a process called mutation. In certain cases, a sequence of mutations in a 
cluster algebra can correspond to iteration of a birational map, so that a discrete dynamical 
system is generated. The reason why cluster algebras have attracted so much attention is that cluster mutations and 
associated discrete dynamical systems or difference equations arise in such a wide variety of contexts, 
including Teichmuller theory \cite{f2,FT}, Poisson geometry \cite{gsv1}, representation theory \cite{difk}, and 
integrable models in  statistical mechanics and quantum field theory \cite{eager,gk,zam}, to name but a few. 

The purpose of this review is to give a brief introduction 
to cluster algebras, and describe certain situations where the 
associated dynamics is completely integrable, in the sense that a discrete version of Liouville's 
theorem in classical mechanics is valid. Furthermore, within the context of cluster algebras, 
we will describe a  way to detect whether a given discrete system is integrable, based on an associated 
tropical dynamical system and its connection to the notion of algebraic entropy. Finally, we describe how discrete 
Painlev\'e equations can arise in the context of cluster algebras.

\section{Cluster algebras: definition and examples} 

A cluster algebra with coefficients, of rank $N$,  is generated by starting from a seed 
$(B, {\bf x}, {\bf y})$ consisting of an \textit{exchange matrix} $B=(b_{ij}) \in \mathrm{Mat}_N(\Z)$, an 
$N$-tuple of \textit{cluster variables} ${\bf x} = (x_1,x_2,\ldots, x_N)$, 
and another $N$-tuple of \textit{coefficients} ${\bf y} = (y_1,y_2,\ldots, y_N)$. The exchange 
matrix is assumed to be skew-symmetrizable, meaning that there is a diagonal matrix $D$, 
consisting of positive integers, such that $DB$ is skew-symmetric. For each integer 
$k\in [1,N]$, there is a  mutation $\mu_k$ which produces a new seed 
$(B', {\bf x}', {\bf y}')=\mu_k (B, {\bf x}, {\bf y})$. The mutation $\mu_k$ consists of
three parts: \textit{matrix mutation},
which is applied to  $B$ to 
produce $B'=(b_{ij}')=\mu_k(B)$, where 
\beq\label{matmut}
b_{ij}' =  \begin{cases}
 -b_{ij} &\text{if}  \,\,i=k \,\, \text{or} \,\, j=k , \\
 b_{ij}+\text{sgn}(b_{ik})[b_{ik}b_{kj}]_+ & \text{otherwise}, 
\end{cases} 
\eeq
with 
$
\text{sgn}(a)$ being $\pm1$ for positive/negative $a\in\R$ and 0 for $a=0$, 
and 
$$ 
[a]_+ =\max (a,0);
$$ 
\textit{coefficient mutation}, defined by ${\bf y}'=(y_j')=\mu_k({\bf y})$ where 
\beq \label{comut} 
y_j'=\begin{cases} 
y_k^{-1} & \text{if} \,\, j=k, \\ 
y_j \left(1+y_k^{-\text{sgn}(b_{jk})} \right)^{-b_{jk}}  & \text{otherwise}; 
\end{cases} 
\eeq 
and \textit{cluster mutation}, given by ${\bf x}'=(x_j')=\mu_k({\bf x})$ with 
the \textit{exchange relation} 
\beq \label{clmut} 
x_k'=
\frac{
y_k\prod_{i=1}^N x_i^ {[b_{ki}]_+} 
+\prod_{i=1}^N x_i^ {[-b_{ki}]_+} 
}
{(1+y_k)x_k} , 
\eeq 
and   
$x_j'= x_j$ for $ j\neq k$.

Given an initial seed, one can apply an arbitrary sequence of mutations, which produces a 
sequence of seeds. This can be visualized by attaching the initial seed to the root of 
an $N$-regular tree $\T_N$ (with $N$ branches attached to each vertex), and then labelling 
the seeds as  $(B_{\bf t}, {\bf x}_{\bf t}, {\bf y}_{\bf t})$ with ``time'' ${\bf t}\in \T_N$.
Note that mutation is an involution, $\mu_k\cdot \mu_k=\text{id}$, but in general two 
successive mutations do not commute, i.e.\  typically $\mu_j\cdot \mu_k\neq \mu_k\cdot \mu_j$ 
for $j\neq k$. Moreover, in general the exponents and coefficients appearing 
in the exchange relation (\ref{clmut}) change at each stage, because the matrix $B$ and the ${\bf y}$ 
variables are altered by each of the previous mutations.  

\begin{definition}\label{ca}  The cluster algebra 
${\cal A}(B, {\bf x}, {\bf y})$ is the algebra over $\C ({\bf y})$ generated by 
the cluster variables produced by all possible sequences of mutations 
applied to  the seed  $(B, {\bf x}, {\bf y})$.
\end{definition}

We will also consider the case of coefficient-free cluster algebras, for which the  ${\bf y}$ 
variables are absent, the seeds are just $(B,{\bf x})$, and the cluster mutation 
is defined by the simpler exchange relation 
\beq \label{freeclmut} 
x_k'=
\frac{
\prod_{i=1}^N x_i^ {[b_{ki}]_+} 
+\prod_{i=1}^N x_i^ {[-b_{ki}]_+} 
}
{x_k} . 
\eeq 

\begin{remark} The original definition of a cluster algebra in \cite{fz1}  involves a more general setting in which 
the coefficients ${\bf y}$ are elements of a semifield $\Ps$,  that is, an abelian multiplicative group together 
with a binary operation $\oplus$ that is commutative, associative and distributive with respect to multiplication. 
In that setting, with the $N$-tuple ${\bf y}\in \Ps^N$, 
the algebra ${\cal A}(B, {\bf x}, {\bf y})$ is defined over $\Z [\Ps ]$, and the addition in the denominator 
of (\ref{clmut}) is given by $\oplus$. The case we consider here corresponds to $\Ps=\Ps_{univ}$, the universal 
semifield, consisting of subtraction-free rational functions in the variables $y_j$, in which case $\oplus$ 
becomes ordinary addition in the field of rational functions $\C ({\bf y})$. However, starting with the more general 
setting, we can also consider the case of the trivial semifield with one element, $\Ps =\{ 1\}$, which yields 
the coefficient-free case (\ref{freeclmut}).  
\end{remark} 


In order to illustrate the above definitions, we now present a number of concrete examples. For the sake of simplicity, 
we concentrate on the coefficient-free case in the rest of this section, and return to the equations with coefficients 
${\bf y}$ at a later stage. 

\begin{example} \label{b2} {\bf The cluster algebra  of type $B_2$:} A particular cluster algebra of rank 
$N=2$ is given by taking the exchange matrix 
\beq\label{b2b} 
B= \left( \begin{array}{cc} 0 & 2 \\ 
-1 & 0 \end{array} \right), 
\eeq 
and the initial cluster ${\bf x}=(x_1,x_2)$, to define a seed $(B,{\bf x})$. 
The matrix $B$ is skew-symmetrizable: the diagonal matrix  $D=\mathrm{diag}(1,2)$ is such that 
$$ 
DB=\left( \begin{array}{cc} 0 & 2 \\ 
-2& 0 \end{array} \right)
$$ 
is skew-symmetric. 
Applying the mutation $\mu_1$ and using the rule (\ref{matmut}) gives a new exchange matrix 
$$ 
B'=\mu_1(B) =   \left( \begin{array}{cc} 0 & -2 \\ 
1 & 0 \end{array} \right)=-B,
$$  
while the coefficient-free exchange relation (\ref{freeclmut}) gives a new cluster 
${\bf x} '=(x_1',x_2)$ with 
$$ 
x_1'=\frac{x_2^2 +1}{x_1}.
$$ 
Since mutation acts as an involution, we have 
$\mu_1 (B',{\bf x}') =(B,{\bf x})$, so  nothing new is obtained by applying $\mu_1$ to this new seed. Thus we 
consider  $\mu_2 (B',{\bf x}')=\mu_2\cdot \mu_1 (B,{\bf x})$ instead, which produces $\mu_2(B')=B$ and $\mu_2({\bf x}') =(x_1',x_2')$, where 
$$ 
x_2'=\frac{x_1'+1}{x_2} = \frac{x_1+x_2^2+1}{x_1x_2}. 
$$ 
Once again, a repeat application of the same mutation $\mu_2$ returns to the previous seed, so instead we consider 
applying $\mu_1$ to obtain $\mu_1\cdot \mu_2\cdot \mu_1(B)=-B$ and 
 $\mu_1\cdot \mu_2\cdot \mu_1({\bf x}) = (x_1'',x_2')$, with 
$$ 
x_1'' = \frac{x_1^2 + 2x_1+x_2^2+1}{x_1x_2^2}.
$$ 
Repeating this sequence of mutations, it is clear that the exchange matrix just changes by an overall sign at each step. 
Perhaps more surprising is the  fact that after obtaining 
 $( \mu_2\cdot \mu_1)^2({\bf x}) = (x_1'',x_2'')$, with 
$$ 
x_2''= \frac{x_1+1}{x_2}, 
$$ the variable $x_1$ reappears in the cluster after 
a further step, i.e.\ 
$\mu_1\cdot ( \mu_2\cdot \mu_1)^2({\bf x}) = (x_1,x_2'')$, and finally 
 $( \mu_2\cdot \mu_1)^3({\bf x}) = (x_1,x_2)={\bf x}$, so that the initial seed  $(B,{\bf x})$ 
is restored after a total of six mutations. Thus the cluster algebra has a finite number of  generators  in this case, 
since there are only 
the six cluster variables $x_1,x_2,x_1',x_2',x_1'',x_2''$. This example is called the cluster 
algebra of type $B_2$, since the initial matrix $B$ is derived from the Cartan matrix of the $B_2$ root system, that is 
$$ C=   \left( \begin{array}{cc} 2 & -2 \\ 
-1 & 2 \end{array} \right), 
$$ 
by replacing the diagonal entries in $C$ with 0, and changing signs of the off-diagonal entries so that 
$b_{ij}$ and $b_{ji}$ have opposite signs for $i\neq j$. 
\end{example} 

There are two significant features of the preceding example, namely the fact that there are only finitely many 
clusters, and the fact that the cluster variables are all Laurent polynomials (polynomials in $x_1$, $x_2$ and their 
reciprocals) with integer coefficients. The first feature is rare: a cluster algebra is said to be of finite type if there are 
only finitely many clusters, and it was shown in \cite{fz2} that all such cluster algebras are generated from seeds 
corresponding to the finite root systems that appear in the Cartan-Killing classification of finite-dimensional 
semisimple Lie algebras. The second feature (the Laurent phenomenon) is ubiquitous \cite{fz}, and follows from 
the following 
result, proved in \cite{fz}. 

\begin{proposition}\label{laur}  All cluster variables in a coefficient-free cluster algebra ${\cal A}(B, {\bf x})$ 
are Laurent polynomials in the variables from the initial cluster, with integer coefficients, i.e.\ they 
are elements of the ring of Laurent polynomials, that is 
$\Z[{\bf x}^{\pm 1}]:=\Z[x_1^{\pm 1}, x_2^{\pm 1},\ldots,x_N^{\pm 1}]$. 
\end{proposition}

\begin{figure} \centering 
\includegraphics[width=5cm,height=5cm,keepaspectratio]{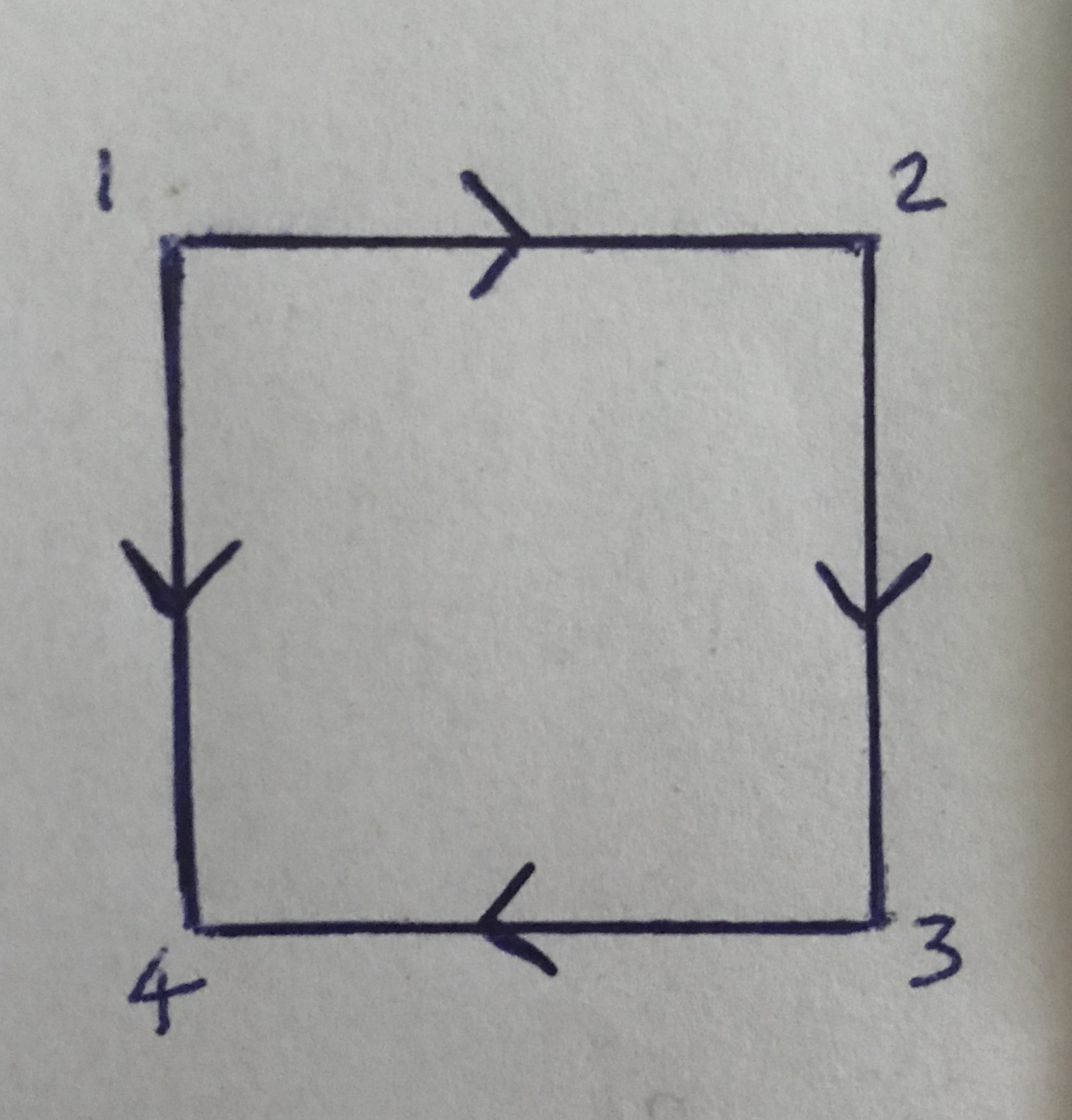}
\caption{\small{ 
The quiver $Q$ corresponding to the exchange matrix (\ref{a13b}).  }}
\label{q1}
\end{figure}

There is an analogous statement in the case that coefficients are included, and in fact it is possible to 
prove the stronger result that all of the coefficients of the cluster variables have positive integer coefficients, so 
they belong to $\Z_{>0}[{\bf x}^{\pm 1}]$ (see \cite{GHKK,ls}, for instance).

\begin{figure} \centering 
\includegraphics[width=5cm,height=5cm,keepaspectratio]{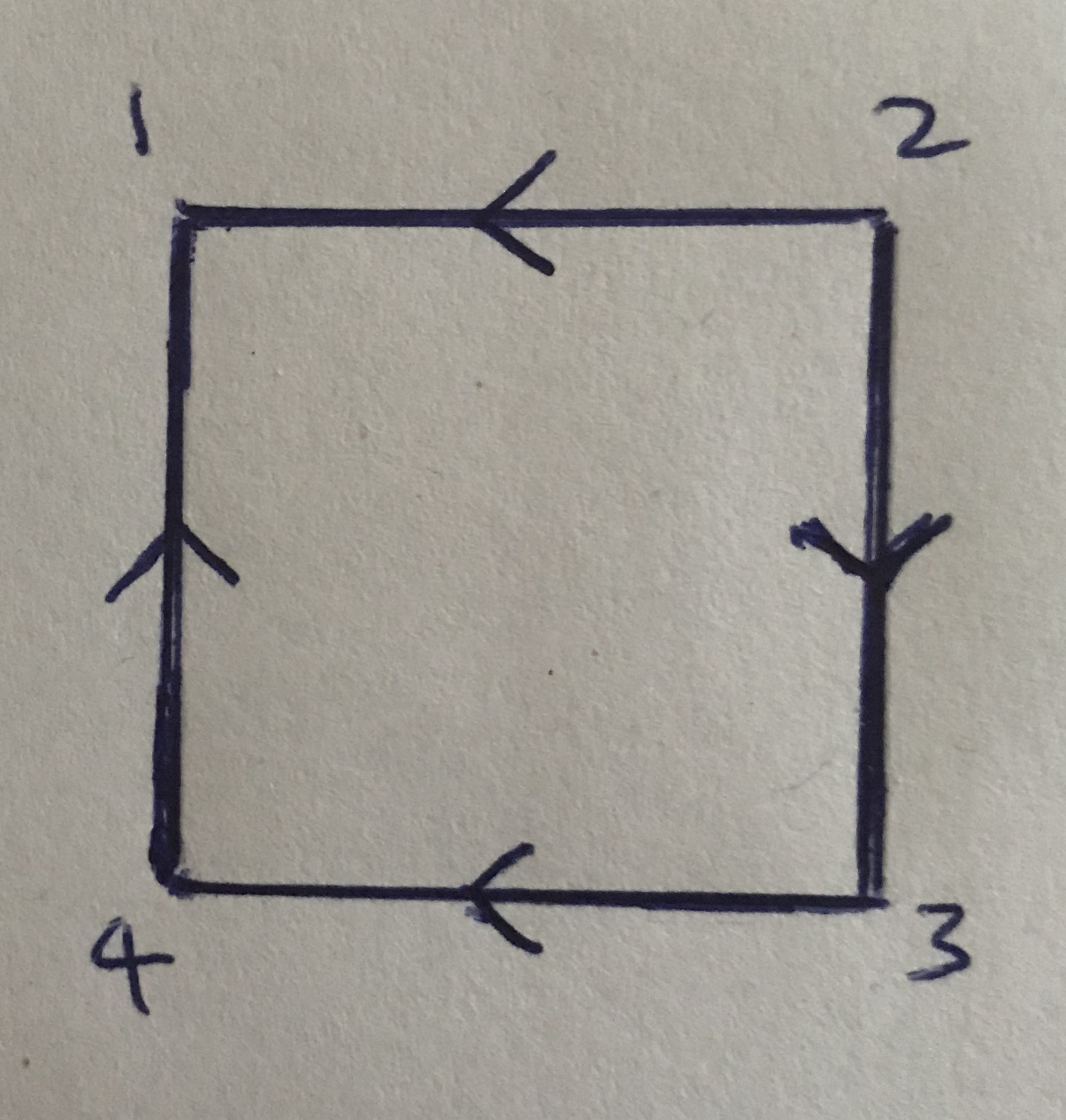}
\caption{\small{ 
The mutated quiver $Q'=\mu_1(Q)$ obtained by applying $\mu_1$ to  (\ref{a13b}).  }}
\label{q2}
\end{figure}

\begin{example} \label{a13} {\bf The cluster algebra  of type $\tilde{A}_{1,3}$:} As an example of rank $N=4$, 
we take the skew-symmetric matrix
\beq\label{a13b} 
B= \left( \begin{array}{cccc} 0 & 1 & 0 & 1 \\ 
-1 & 0 & 1 & 0 \\ 
0 & -1 & 0 & 1 \\ 
-1 & 0 & -1 & 0 
\end{array} \right), 
\eeq 
which is obtained from the Cartan matrix of the  affine root system $A_3^{(1)}$ \cite{kac}, namely 
$$ 
C= \left( \begin{array}{cccc} 2 & -1 & 0 & -1 \\ 
-1 & 2 & -1 & 0 \\ 
0 & -1 & 2 & -1 \\ 
-1 & 0 & -1 & 2 
\end{array} \right),
$$ 
by replacing each of the diagonal entries of $C$ with 0, and making a suitable adjustment of signs 
for the off-diagonal entries, such that $b_{ij}=-b_{ji}$. Since $B$ is a skew-symmetric integer matrix, it can be associated with a 
quiver $Q$ without 1- or 2-cycles, that is, a directed graph specified by the rule that 
$b_{ij}$  is equal to the number of arrows $i\rightarrow j$ if it is non-negative, and minus the number of arrows $j\rightarrow i$ otherwise (see Fig.\ref{q1}). 
If the mutation $\mu_1$ is applied, then the new exchange matrix is 
$$ 
B'=\mu_1(B)= \left( \begin{array}{cccc} 0 & -1 & 0 & -1 \\ 
1 & 0 & 1 & 0 \\ 
0 & -1 & 0 & 1 \\ 
1 & 0 & -1 & 0 
\end{array} \right), 
$$ 
which corresponds to a new quiver $Q'$ obtained by a cyclic permutation of the vertices of the original $Q$ 
 (see Fig.\ref{q2}), 
while the initial cluster ${\bf x}=(x_1,x_2,x_3,x_4)$ is mutated to 
$\mu_1({\bf x})=(x_1',x_2,x_3,x_4)$, where $x_1'$ is defined by the relation 
$$ 
x_1 x_1' = x_2 x_4 +1. 
$$ 
Rather than trying to describe the effect of every possible choice of mutation, we 
consider what happens when $\mu_1$ is followed by $\mu_2$, and once more 
observe that, at the level of the associated quiver, this just corresponds to 
applying the same cyclic permutation as before to the vertex labels 1,2,3,4.   
The new cluster obtained from this is $\mu_2\cdot \mu_1 ({\bf x})=(x_1',x_2',x_3,x_4)$, with $x_2'$ 
defined by 
$$ 
x_2 x_2'=x_3 x_1'+1,
$$ 
and if $\mu_3$ is applied next, then $\mu_3\cdot \mu_2\cdot \mu_1 ({\bf x})=(x_1',x_2',x_3',x_4)$, with
$$ 
x_3 x_3'=x_4 x_2'+1.
$$
Continuing in this way, it is not hard to see that the composition $\mu_4\cdot\mu_3\cdot\mu_2\cdot \mu_1$ 
takes the original $B$ to itself, and applying this sequence of mutations repeatedly in the same order generates a 
new cluster variable at each step, 
with the sequence of cluster variables satisfying the nonlinear recurrence relation 
\beq\label{a13rec} 
x_nx_{n+4}=x_{n+1}x_{n+3}+1 
\eeq 
(where we have made the identification $x_1'=x_5$, $x_2'=x_6$, and so on).  
Regardless of other possible choices of mutations, this particular sequence of mutations alone generates 
an infinite set of distinct cluster variables, as can be seen by fixing some numerical values for the initial cluster.
In fact, as was noted in \cite{hone}, for any orbit of (\ref{a13rec}) there is a constant $K$ such that the 
iterates satisfy the 
linear recurrence 
\beq\label{linrec} 
x_{n+6}+x_n = Kx_{n+3}.
\eeq  
Upon fixing $(x_1,x_2,x_3,x_4)=(1,1,1,1)$, the nonlinear recurrence generates the integer sequence 
$$ 
1,1,1,1,2,3,4,9,14,19,43,76,\ldots , 
$$  
which also satisfies the linear recurrence (\ref{linrec}) with $K=5$; so the terms grow 
exponentially with $n$, and the integers $x_n$  are distinct for $n\geq 4$. 
This is called the $\tilde{A}_{1,3}$ cluster algebra, because the corresponding quiver is 
 an orientation of the edges of an   
affine Dynkin diagram of type $A$ with one anticlockwise  arrow and three clockwise
 arrows.
\end{example}

The skew-symmetry of $B$ is preserved under matrix mutation, and for any skew-symmetric integer matrix 
there is an equivalent operation of \textit{quiver mutation} which acts on the associated quiver $Q$: to 
obtain the mutated quiver $\mu_k(Q)$ one should (i) add $pq$ arrows $i\overset{pq}{\longrightarrow}j$ 
whenever $Q$ has a path of length two passing through vertex $k$ with $p$ arrows $i\overset{p}{\longrightarrow} k$
and $q$ arrows 
$k\overset{q}{\longrightarrow} j$; (ii) reverse all arrows in $Q$ that go in/out of vertex $k$; (iii) delete 
any 2-cycles created in the first step.

Unlike the $B_2$ cluster algebra, 
the above example is not of finite type, because there are infinitely many clusters. However, it turns out that it is of 
finite mutation type, in the sense that 
there are only a finite number of exchange matrices produced under mutation from the initial $B$. 
Cluster algebras of finite mutation type have also been classified \cite{fst,fst2}: as well as those of finite type, they 
include cluster algebras associated with triangulated surfaces \cite{f2,FT}, 
cluster algebras of rank 2, 
plus a finite number of exceptional cases.  

\begin{figure} \centering 
\includegraphics[width=5cm,height=5cm,keepaspectratio]{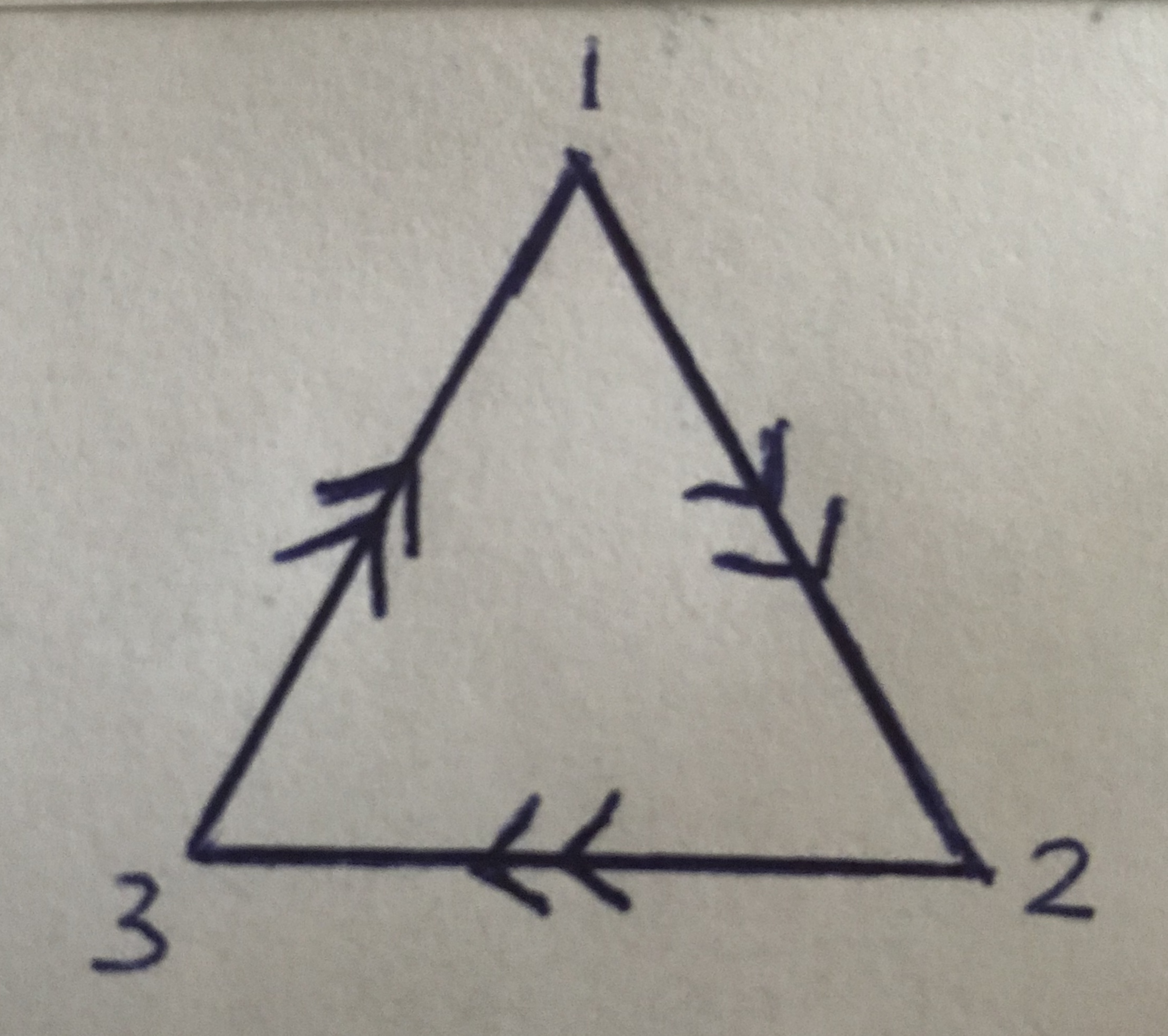}
\caption{\small{ 
The quiver corresponding to the exchange matrix  (\ref{markoffb}).  }}
\label{q3}
\end{figure}

\begin{example}\label{markoff}  {\bf Cluster algebra related to Markoff's equation:} 
For $N=3$, consider the exchange matrix 
\beq \label{markoffb}
B=\left(\begin{array}{ccc} 0 & 2 & -2 \\ 
-2 & 0 & 2 \\ 
2 & -2 & 0 \end{array} \right),
\eeq 
which is associated with the quiver in Fig.\ref{q3}. 
After any sequence of matrix mutations, one can obtain only $B$ or $-B$, so this is another example 
of finite mutation type: it is connected to the moduli space of once-punctured tori, and the Markoff 
equation 
\beq \label{meq} 
x^2 + y^2 + z^2 = 3xyz
\eeq 
which arises in that context as well as in Diophantine approximation theory \cite{cassels,cohn}. Upon 
applying $\mu_1$ to the initial cluster $(x_1,x_2,x_3)$, the result is 
$(x_1',x_2,x_3)$ with 
$$
x_1x_1'=x_2^2 +x_3^2,
$$ 
and a subsequent application of $\mu_2$ yields $(x_1',x_2',x_3)$, where  
$$ 
x_2 x_2'=x_3^2 + x_1'^2.
$$ 
Repeated application of the  mutations $\mu_3\cdot\mu_2\cdot \mu_1$ in that order produces a new cluster 
variable at each step, and upon identifying $x_4=x_1'$, $x_5=x_2'$, and so on, the sequence of cluster variables 
$(x_n)$ is generated by a  recurrence of third order, namely 
\beq\label{mrec} 
x_nx_{n+3} = x_{n+1}^2 + x_{n+2}^2 .
\eeq  
It can also be shown that on each orbit of (\ref{mrec}) there is a constant $K$ such that the 
nonlinear relation 
$$ 
x_{n+3} + x_n = Kx_{n+1}x_{n+2} 
$$ 
holds for all $n$, and by using the latter to eliminate $x_{n+3}$ it follows that 
\beq\label{km} 
K = \frac{x_n^2+x_{n+1}^2 +x_{n+2}^2}{x_nx_{n+1}x_{n+2}} 
\eeq
is an invariant for (\ref{mrec}), independent of $n$. In particular, 
taking the initial values to be $(1,1,1)$ gives $K=3$, and each adjacent triple 
$(x,y,z)=(x_n,x_{n+1},x_{n+2})$ in the resulting sequence   
\beq\label{mseq}
1,1,1,2,5,29,433,37666,48928105,\ldots 
\eeq 
is an integer solution of Markoff's equation (\ref{meq}). 
The terms of this sequence have double exponential growth: $\log x_n$ grows exponentially with $n$.
\end{example} 

The next example is generic, in the sense that there are both infinitely many clusters and infinitely many exchange matrices.  

\begin{figure} \centering 
\includegraphics[width=5cm,height=5cm,keepaspectratio]{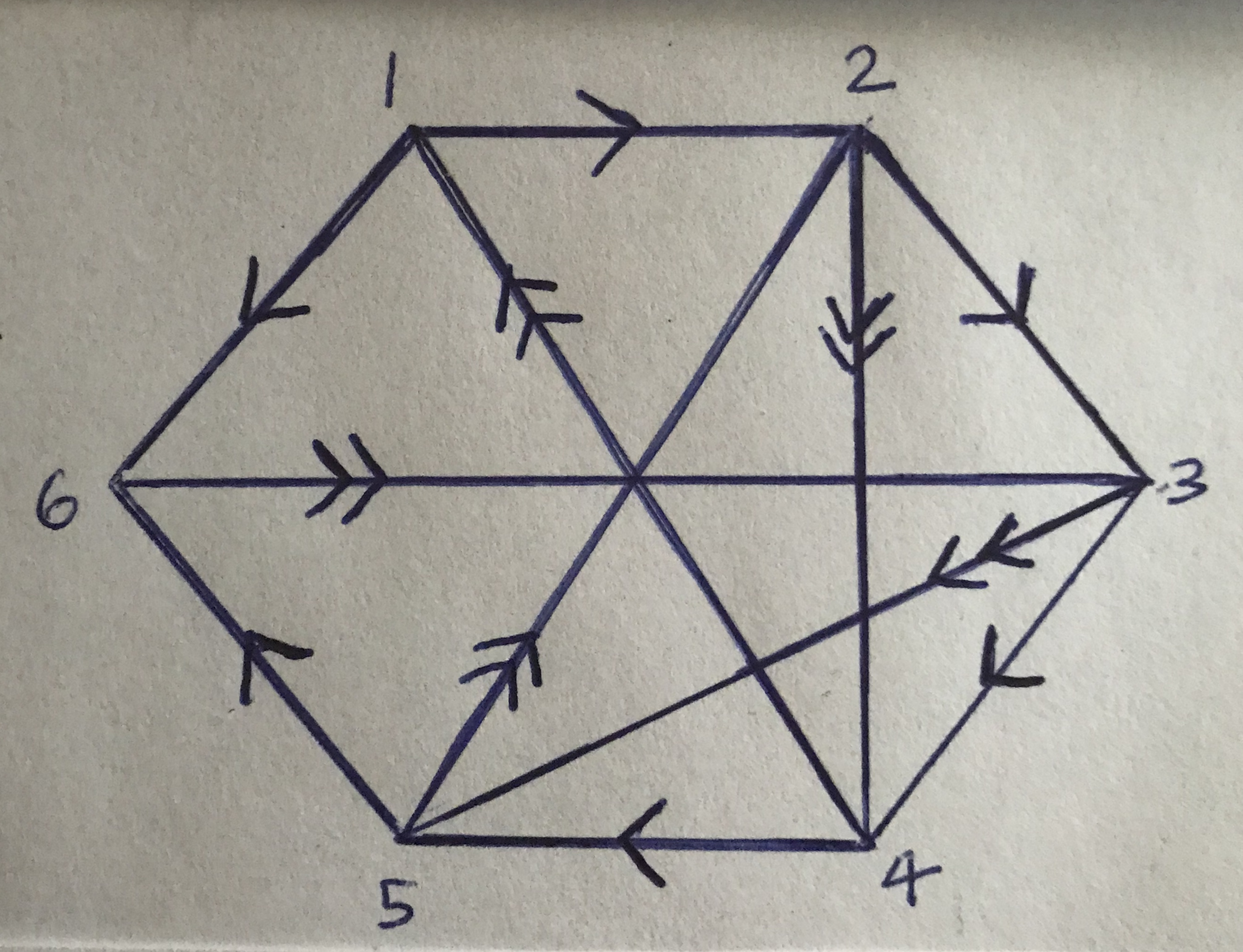}
\caption{\small{ 
The Somos-6 quiver corresponding to the exchange matrix  (\ref{s6b}).  }}
\label{q4}
\end{figure}
\begin{example}\label{s6} {\bf A Somos-6 recurrence:} A  sequence  
that is generated   by a quadratic recurrence relation of the form 
$$ 
x_{n}x_{n+k} = \sum_{j=1}^{\left \lfloor{k/2}\right \rfloor} \alpha_j \, x_{n+j}x_{n+k-j}, 
$$ 
where $\alpha_j$ are coefficients, 
is called a Somos-$k$ sequence (see \cite{fz,gale,honetams,svdp,rob}). 
A certain class of Somos-6 sequences can be generated by 
starting from the exchange matrix 
\beq\label{s6b} 
B=\left(\begin{array}{cccccc} 
0 & 1 & 0 & -2 & 0 & 1 \\ 
-1 & 0 & 1 & 2 & -2 & 0 \\ 
0 & -1 & 0 & 1 & 2 & -2 \\ 
2 & -2 & -1 & 0 & 1 & 0 \\ 
0 & 2 & -2 & -1 & 0 & 1 \\ 
-1 & 0 & 2 & 0 & -1 & 0 \end{array}\right), 
\eeq 
which corresponds to the quiver in Fig.\ref{q4}.
Upon applying cyclic sequences of mutations ordered as 
$\mu_6\cdot\mu_5\cdot\mu_4\cdot\mu_3\cdot\mu_2\cdot\mu_1$, a sequence of cluster variables $(x_n)$ 
is produced which 
satisfies the particular  Somos-6 recurrence 
\beq\label{s6rec}
x_{n}x_{n+6} =  x_{n+1}x_{n+5} + x_{n+3}^2. 
\eeq 
If six 1s are chosen as initial values, then an integer Somos-6 sequence beginning with 
$$ 
1,1,1,1,1,1,2,3,4,8,17,50, 107, 239,\ldots 
$$ 
is produced. For this sequence, $\log x_n$ grows like $n^2$. However, applying successive mutations other than these cyclic ones  generally causes the magnitude  of the entries of the exchange matrices to grow - for instance, 
$$ 
\mu_5\cdot\mu_4\cdot\mu_2 (B)=\left( \begin{array}{cccccc} 
0 & -1 & 1 & 0 & 0 & 1 \\ 
1 & 0 & -3 & -10 & 4 & 0 \\ 
-1 & 3 &  0 & -1 & 0 & -2 \\ 
0 & 10 & 1 & 0 & -3 & 3 \\ 
0 & -4 & 0 & 3 & 0 & -1 \\ 
-1 & 0 & 2 & -3 & 1 & 0 
\end{array}\right); 
$$ 
 and 
(e.g.\ starting with the initial seed evaluated as $(1,1,1,1,1,1)$ as before) typically this results in the values of cluster variables 
showing  double 
exponential growth with the number of steps. 
\end{example}

\section{Cluster algebras with periodicity} 

The exchange relation (\ref{clmut}) can be regarded as a birational map in $\C^N$. Alternatively, 
${\bf x}\in (\C^*)^N$ can be viewed as coordinates in a toric chart for some algebraic variety, and a mutation 
${\bf x} \mapsto \mu_k({\bf x})={\bf x}'$ as a change of coordinates to another chart. 
The latter point of view is passive, in the sense that there is some fixed variety and mutation  just selects  
different choices of coordinate charts. Instead of this, we would like to take an active view,  regarding 
each mutation as an iteration in a discrete dynamical system. However, there is a problem with this, because 
a general sequence of mutations is specified by a ``time'' ${\bf t}$ belonging to the tree $\T_N$, 
and   (except for the case of rank $N=2$) this cannot naturally be identified with a discrete time belonging 
to the set of  integers $\Z$.  Furthermore, there is the additional problem that matrix mutation, as in 
(\ref{matmut}), typically changes the exponents appearing in the exchange relation, so that in general 
it is not possible to interpret successive mutations as iterations of the same map. 

Despite the above comments, it turns out that the most interesting cluster algebras appearing ``in nature'' 
have special symmetries, in the sense that they display periodic behaviour with respect to at least 
some subset of the possible mutations. In fact, all of the examples in the previous section are of this kind. 
Here we consider a notion of periodicity that was introduced by Fordy and Marsh \cite{fm} in the context of 
skew-symmetric exchange matrices $B$, which correspond to quivers. 

\begin{definition}\label{per} 
An exchange matrix $B$ is said to be {\it cluster mutation-periodic} with period $m$ if (for a suitable labelling 
of indices) $\mu_m\cdot \mu_{m-1} \cdot \ldots \cdot \mu_1 (B) = \rho^m (B)$, where $\rho$ 
is the cyclic permutation $$\rho:\,(1,2,3,\ldots,N)\mapsto (N,1,2,\ldots,N-1).$$   
\end{definition}  

In the context of quiver mutation, the case of cluster mutation-periodicity with period $m=1$ means that 
the action of mutation $\mu_1$ on $Q$ is the same as the action of $\rho$, which is such that the number 
of arrows $i \rightarrow j$ in $Q$ is the same as the number of arrows $\rho^{-1}(i)\rightarrow \rho^{-1}(j)$ in 
$\rho(Q)$. This means that the cluster map $\varphi=\rho^{-1}\cdot \mu_1$ acts as the identity on $Q$ (or equivalently, 
on $B$), but in general ${\bf x}\mapsto \varphi ({\bf x})$ has a non-trivial action on the cluster. Mutation-periodicity 
with period 1 implies that iterating this map is equivalent to iterating a single recurrence relation. 

\begin{example} 
Although it is not skew-symmetric, the exchange matrix (\ref{b2b}) in Example \ref{b2} is cluster mutation-periodic with period 2, since $\mu_2\cdot\mu_1(B) = \rho^2(B) = B$ where $\rho$ is the switch $1\leftrightarrow 2$. 
Defining the cluster map to be $\varphi=\rho^{-2}\cdot \mu_2\cdot\mu_1$, the action 
of $\varphi$ is 
periodic with period 3 for any choice of initial cluster, i.e.\ $\varphi^3=$id. 
\end{example} 

\begin{example} The exchange matrix (\ref{a13b}) in Example \ref{a13} is cluster mutation-periodic with period 1. 
The cluster map $\varphi=\rho^{-1}\cdot \mu_1$ is given by 
\beq\label{a13map} 
\varphi: \,\, (x_1,x_2,x_3,x_4)\mapsto \left( x_2,x_3,x_4, \frac{x_2x_4+1}{x_1} \right) ,
\eeq 
whose iterates are equivalent to those of the nonlinear recurrence (\ref{a13rec}). 
\end{example} 

\begin{example} The exchange matrix (\ref{markoffb}) in Example \ref{markoff} is cluster mutation-periodic with period 2. 
The cluster map $\varphi=\rho^{-2}\cdot \mu_2\cdot \mu_1$ is given by 
\beq\label{mmap} 
\varphi:\,\,(x_1,x_2,x_3)\mapsto \left( x_3,x_1', \frac{x_3^2+x_1'^2}{x_2}\right) ,\quad \mathrm{with} \quad 
x_1'=\frac{x_2^2+x_3^2}{x_1}.
\eeq 
Each  iteration of (\ref{mmap}) is equivalent to two iterations of the nonlinear recurrence (\ref{mrec}). 
This period 2 example is exceptional because, in general, cluster mutation-periodicity with period $m>2$ does not give rise to a single recurrence relation 
(see \cite{fm} for more examples).
\end{example} 

\begin{example} The exchange matrix (\ref{s6b}) in Example \ref{s6} is cluster mutation-periodic with period 1. 
The cluster map is given by 
\beq\label{s6map} 
\varphi : \, 
(x_1,x_2,x_3,x_4,x_5,x_6)\mapsto \left( x_2,x_3,x_4,x_5,x_6, \frac{x_2x_6+x_4^2}{x_1}\right) ,
\eeq 
whose iterates are equivalent to those of the nonlinear recurrence (\ref{s6rec}). 
\end{example} 

\begin{remark} There is a more general notion of periodicity, due to Nakanishi \cite{nak}, which extends Definition \ref{per}. This yields broad
generalizations of Zamolodchikov's Y-systems  \cite{zam}, a set of functional relations, arising from the 
thermodynamic Bethe ansatz for certain integrable quantum field theories, that were the prototype 
for the coefficient mutation (\ref{comut}) in a cluster algebra. We shall introduce examples of 
generalized Y-systems in the sequel.  
\end{remark}

Fordy and Marsh gave a complete classification of period 1 quivers. Their result can be paraphrased as follows. 

\begin{theorem} \label{peri}
Let $(a_1,\ldots, a_{N-1})$ an $(N-1)$-tuple of integers that is palindromic, i.e.\ $a_j=a_{N-j}$ for all $j\in [1,N-1]$.  
Then the skew-symmetric exchange matrix $B=(b_{ij})$ with entries specified by 
$$ 
b_{1,j+1}=a_j
\quad
and 
\quad 
b_{i+1,j+1}=b_{ij} + a_i[-a_j]_+ - a_j[-a_i]_+,
$$ 
for all $i,j\in[1,N-1]$, is cluster mutation-periodic with period 1, and every period 1 skew-symmetric $B$ arises 
in this way. 
\end{theorem} 

The above result says that a period 1 skew-symmetric $B$ matrix is completely determined by the entries in its first row 
(or equivalently, its first column), and these form a palindrome after removing $b_{11}$. The 
entries $a_j$ in the palindrome are precisely the exponents that appear in the exchange relation defining the cluster 
map $\varphi$, whose iterates are equivalent to those of the nonlinear recurrence relation 
\beq\label{perrec} 
x_n x_{n+N} = \prod_{j:\, a_j>0} x_{n+j}^{a_j} + \prod_{j:\, a_j<0} x_{n+j}^{-a_j}.
\eeq 
Thus (\ref{perrec}) corresponds to a special sequence of mutations in a particular subclass of cluster algebras. Such a  nonlinear recurrence is an example of a generalized T-system, in the terminology of \cite{nak}. 
 
Next we would like to turn to the question of which recurrences of this special type correspond to discrete 
integrable systems. We begin our approach to this question in the next section, by considering the notion of algebraic entropy, which 
gives a measure of the growth of iterates in a discrete dynamical system defined by iteration of rational functions. 

\section{Algebraic entropy and tropical dynamics} 

There are  various different ways of quantifying the growth, or complexity, of a discrete dynamical system (see \cite{aba}, 
for instance). In the context of discrete integrability of birational maps, Bellon and Viallet introduced the concept of 
algebraic entropy, and proposed that zero algebraic entropy should be a criterion for integrability \cite{bellon}. For a birational map $\varphi$, one can calculate the degree $d_n = \deg \varphi^n$, given 
by the maximum of the degrees of the components of the map $\varphi^n$, and then the algebraic 
entropy is defined to be 
$$ 
{\cal E}:=\lim_{n\to\infty} \frac{\log d_n}{n}. 
$$ 
Typically, the degree $d_n$ grows exponentially with $n$, so ${\cal E}>0$, but in rare cases there can be subexponential growth, leading to vanishing entropy. In the case of birational maps in two dimensions, 
the types of degree growth have been fully classified \cite{df}, and there are only four possibilities: bounded degrees, linear growth, quadratic growth, or exponential growth; the first three cases, with zero entropy, coincide with 
the existence of invariant foliations. Thus, at least for maps of the plane, the requirement of zero entropy identifies 
symplectic maps that are integrable in the sense that they satisfy the 
conditions needed for a discrete analogue of the Liouville-Arnold theorem to hold   \cite{bruschi,maeda,veselov}.

 Measuring the degree growth and seeking maps with zero algebraic entropy is a useful tool for identifying 
discrete integrable systems. (For another approach, based on the growth of heights in orbits defined 
over $\Q$ or a number field, see \cite{halburd}.) Once such a map has been identified, it leaves open the 
question of Liouville integrability; this is discussed in the next section. For now, we concentrate on 
the case of maps arising from cluster algebras, and consider algebraic entropy in that setting. 

The advantage of working with cluster maps is that, due to the Laurent property, it is sufficient to consider the growth 
of degrees of the  denominators of the  cluster variables in order to determine the algebraic entropy. 
In particular, in the period 1 case, by 
Proposition \ref{laur} every iterate of (\ref{perrec}) can be written in the form 
\beq\label{lpoly} 
x_n = 
\frac{ \mathrm{P}_n (\bx) } 
{ \bx^{\bd_n} } , 
\eeq 
where the polynomial $\mathrm{P}_n$ is not divisible by any of the $x_j$ from the initial cluster, 
and the monomial 
$ \bx^{\bd_n} =\prod_{j=1}^N x_j^{ d_n^{(j)} }$ 
is specified by 
the integer vector $${\bf d}_n=(d_n^{(1)},\ldots,d_n^{(N)})^T,$$ 
known as a d-vector. From the fact that cluster variables are subtraction-free rational expressions in 
${\bf x}$ (or, a fortiori, from the fact that these Laurent polynomials are now known to have positive integer coefficients 
\cite{GHKK}), it follows that the d-vectors in a cluster algebra satisfy the 
max-plus tropical analogue of the exchange relations for the corresponding cluster variables \cite{f2,f1}, 
where the latter is obtained from (\ref{freeclmut}) 
by replacing each addition with max, and each multiplication with addition. 
In the case of (\ref{perrec}), this implies the following result. 

\begin{proposition} 
If $x_n$ given by (\ref{lpoly}) satisfies (\ref{perrec}), then the sequence of vectors ${\bf d}_n$ satisfies the tropical 
recurrence relation 
\beq\label{trop}
{\bf d}_n + {\bf d}_{n+N} = \max \left( \sum_{j:\, a_j>0}a_j {\bf d}_{n+j}, - \sum_{j:\, a_j<0}a_j {\bf d}_{n+j} \right).   
\eeq   
\end{proposition} 
Note that  the equality in (\ref{trop}) holds  componentwise. For a detailed proof of this result, see \cite{FH}.  

In the period 1 situation, the problem of determining the evolution of d-vectors can be simplified further, upon noting that 
the first component of ${\bf d}_n$ has the initial values 
\beq\label{dinits} 
d_1^{(1)} = -1, \qquad 
d_j^{(1)} =0 \quad \mathrm{for} \quad 2\leq j \leq N, 
\eeq 
while each of the other components $d_n^{(k)}$ for $k\in [2,N]$ has the same set of initial values but shifted by $k-1$ steps (so, for instance, $d_1^{(2)}=0$, $d_2^{(2)}=-1$ and then $d_j^{(2)}=0$ for $3\leq j \leq N+1$, since the 
first division by the variable $x_2$ appears in $x_{N+2}$, etc.). The total degree of the 
monomial $\bx^{\bd_n}$ is  the sum 
$\sum_{k=1}^N d_n^{(k)}$ of the components of the 
d-vector, and if the components are all non-negative then this coincides with the degree of the denominator 
of the rational function (\ref{lpoly}). Unless there is periodicity of d-vectors, corresponding to  degrees remaining bounded
(which can only happen in finite type cases like Example \ref{b2}), then all these  components are positive for large 
enough $n$. Moreover, it is not hard to see that the growth of the degree of the numerators $\mathrm{P}_n$ appearing in the Laurent polynomials (\ref{lpoly}) is controlled by that of the denominators. Thus, to   
determine the growth of degrees of Laurent polynomials generated by (\ref{perrec}), it is 
sufficient to consider the solution of the scalar version of (\ref{trop}), with initial data given 
by (\ref{dinits}), and the growth of this determines the algebraic entropy. 

\begin{example}\label{a13g} 
For the recurrence (\ref{a13rec}) in Example \ref{a13}, the tropical equation for determining the 
degrees of d-vectors is given in scalar form by 
\beq \label{tropa13} 
d_n + d_{n+4} = \max (d_{n+1}+d_{n+3},0). 
\eeq 
If we take  initial values $d_1=-1$, $d_2=d_3=d_4=0$, corresponding to (\ref{dinits}), 
then by induction it follows that 
$d_n\geq 0$ for all $n\geq 2$, so that the max on the right-hand side of (\ref{tropa13}) can be replaced by   its 
first entry, to yield the linear recurrence 
$$ 
d_n + d_{n+4} = d_{n+1}+d_{n+3} \quad \mathrm{for} \quad n\geq 1.
$$
The characteristic polynomial of the latter factorizes as $(\la -1)^2 (\la^2 +\la +1)=0$, leading to the 
solution 
$$ 
d_n = n/3 -1 +  \frac{\ri}{3\sqrt{3}} (\varepsilon^n -\varepsilon^{-n}), \quad \varepsilon=(-1+\ri\sqrt{3})/2.
$$ 
Thus we have a   sequence 
that 
grows linearly with $n$, 
beginning with   
$$ 
-1,0,0,0,1,1,1,2, 2, 2,3,3,3,4,4,4,\ldots ,
$$  
where each positive integer appears three times in succession, 
which corresponds to the degree of the denominator of $x_n$ in each of the variables $x_1,x_2,x_3,x_4$ separately. Clearly the total degree of the denominator also grows linearly, and 
the algebraic entropy is 
$\lim_{n\to\infty} (\log d_n)/n =0$ in this case.  
\end{example} 

\begin{example} \label{mg}
The exchange matrix (\ref{markoffb}) in Example \ref{markoff} is period 2 rather than period 1, 
but we can still calculate the growth of d-vectors in the recurrence (\ref{mrec}) by taking its tropical version, 
namely 
$$ 
d_n + d_{n+3}=2\,\max(d_{n+1},d_{n+2}), 
$$ 
and choosing the initial values $d_1=-1$, $d_2=d_3=0$, which produces a sequence beginning 
$$ 
-1,0,0,1,2,4,7,12,20,33,54,88,\ldots .
$$ 
By induction one can show that $d_{n+2}\geq d_{n+1}$ for $n\geq 0$, so in fact 
the linear recurrence 
$$ 
d_{n+3}+d_n = 2d_{n+2} 
$$ 
holds for this sequence, with characteristic equation $(\la -1)(\la^2-\la-1)=0$, and it turns out that the differences 
$$ 
F_n = d_{n+3}-d_{n+2} 
$$ 
are just the Fibonacci numbers. Hence there is a constant C$\,>0$ such that 
$$d_n\sim \mathrm{C}\,\left(\frac{1+\sqrt{5}}{2}\right)^n, 
$$ 
and the algebraic entropy ${\cal E}=\log((1+\sqrt{5})/2)$, which is the same as the limit 
$\lim_{n\to\infty} (\log\log x_n)/n$ for the sequence (\ref{mseq}) - see \cite{dioph}. 
\end{example} 

\begin{example}\label{s6g}
For the period 1 exchange matrix (\ref{s6b}) in Example \ref{s6}, we consider the recurrence 
\beq\label{s6t} 
d_n + d_{n+6}=\max(d_{n+1}+d_{n+5}, 2d_{n+3}), 
\eeq 
which is the max-plus  analogue of (\ref{s6rec}), and take initial data 
\beq\label{ds6} 
d_1=-1, \quad d_2=d_3=\cdots =d_6 = 0, 
\eeq 
which generates a degree sequence beginning 
\beq\label{s6deg} 
-1,0,0,0,0,0,1,1,1,2,2,3,3,3,5,5,6,7,7,9,9,10,12,12,14,15,16,\ldots .
\eeq  
In order to simplify the analysis of (\ref{s6t}), we observe that the combination 
\beq\label{utrop} 
U_n=d_{n+2}-2d_{n+1}+d_n 
\eeq 
satisfies a recurrence of fourth order, namely 
\beq\label{usystrop} 
U_{n+4}+2U_{n+3}+3U_{n+2}+2U_{n+1}+U_n=\max(U_{n+3}+2U_{n+2}+U_{n+1},0). 
\eeq 
(The origin of the substitution (\ref{utrop}) will be explained in the next section.) 
The values in (\ref{ds6}) correspond to the initial conditions 
$$ 
U_1=-1, \quad U_2=U_3=U_4=0
$$ 
for (\ref{usystrop}), which generate a sequence $(U_n)$ beginning with 
$$ \begin{array}{l}
-1,0,0,0,1,-1,0,1,-1,1,-1,0,2,-2,1,0,-1,2,-2,1,1,-2, \\ 
2,-1,0,1,-2,\ldots , 
\end{array}  
$$ 
and further calculation with a computer shows that this sequence does not repeat for the first 40 steps, 
but then $U_{42}=-1$ and 
$U_{43}=U_{44}=U_{45}=0$, so it is periodic with period 41. Thus in terms of the shift operator $\SH$,  
which sends $n\to n+1$, 
$$ 
(\SH^{41}-1) U_n = (\SH^{41}-1)(d_{n+2}-2d_{n+1}+d_n )= (\SH^{41}-1)(\SH-1)^2 d_n=0, 
$$ 
which is a linear recurrence of order 43 satisfied by the degree sequence (\ref{s6deg}). Clearly 
the characteristic polynomial of the latter has $\la=1$ as a triple root, and all other characteristic roots have 
modulus 1. Therefore, for some constant $\mathrm{C}'>0$,  
$$
d_n\sim \mathrm{C}' n^2
$$ 
as $n\to\infty$, which implies that 
(\ref{s6rec}) has algebraic entropy ${\cal E}=0$. 
\end{example} 

The preceding examples indicate that we should regard (\ref{a13rec}) and (\ref{s6rec}) as being integrable in some sense, 
and (\ref{mrec}) as non-integrable. According to the relation (\ref{linrec}), we know that 
(\ref{a13rec}) has at least one conserved quantity $K$; and it turns out to  have three independent  
conserved quantities \cite{hone}. The recurrence (\ref{mrec}) also has a conserved quantity, given by  (\ref{km}), but it is 
possible to show that it can have no other algebraic conserved quantitites, independent of this one. In the next 
section we will derive two independent conserved quantities for (\ref{s6rec}), and we will discuss the 
interpretation of all these examples from the viewpoint of Liouville integrability.   

In \cite{FH}, a detailed analysis of the behaviour of the tropical  recurrences (\ref{trop}) led to the conjecture that 
the algebraic entropy of (\ref{perrec}) should be positive if and only if the following condition holds: 
\beq\label{econ} 
\max\left(\sum_{j=1}^{N-1} [a_j]_+, -\sum_{j=1}^{N-1} [-a_j]_+\right)\geq 3.
\eeq 
In other words, in order for the cluster map defined by (\ref{perrec}) to have a zero entropy, the degree of nonlinearity 
cannot be too large.  The analysis of algebraic entropy for other types of cluster maps has been carried out 
more recently using methods based on Newton polytopes \cite{gp}, and using the same methods it is also possible to prove the above conjecture\footnote{P. Galashin, private communication, 2017}.  By enumerating the 
possible choices of exponents that lie below the bound (\ref{econ}), this leads to a complete proof 
of a classification result for nonlinear recurrences of the form (\ref{perrec}), as stated in \cite{FH}. 

\begin{theorem}\label{main} 
A cluster map $\varphi$ given by a recurrence (\ref{perrec}) has algebraic entropy ${\cal E}=0$ 
if and only if it belongs to one of the following four families: 

\noindent 
(i) For even $N=2m$, recurrences of the form 
\beq\label{i} 
x_nx_{n+2m}=x_{n+m}+1. 
\eeq 

\noindent 
(ii) For $N\geq 2$ and $1\leq q\leq \left \lfloor{N/2}\right \rfloor$, recurrences of the form 
\beq\label{ii} 
x_nx_{n+N}=x_{n+q}x_{n+N-q}+1. 
\eeq

\noindent 
(iii) For even $N=2m$ and $1\leq q\leq m-1$, recurrences of the form 
\beq\label{iii} 
x_nx_{n+2m}=x_{n+q}x_{n+2m-q}+ x_{n+m}. 
\eeq 

\noindent 
(iv) For $N\geq 2$ and $1\leq p<q\leq \left \lfloor{N/2}\right \rfloor$, recurrences of the form  
\beq\label{iv} 
x_nx_{n+N}=x_{n+p}x_{n+N-p} + x_{n+q}x_{n+N-q}. 
\eeq 
\end{theorem}  

Case (i) is somewhat trivial: the recurrence (\ref{i})  is equivalent to taking $m$ copies of the Lyness 5-cycle 
$$ 
x_nx_{n+2}=x_{n+1}+1, 
$$ 
for which every orbit has period 5, corresponding to the cluster algebra of finite type 
associated with the root system $A_2$; so in this case the dynamics is purely periodic and there is no degree growth. 
Both case (ii), which corresponds to affine quivers of type $\tilde{A}_{q,N-q}$, and case (iii) display linear degree growth, 
similar to Example \ref{a13g}. Case (iv) consists of Somos-$N$ recurrences, which display quadratic degree growth \cite{mase2}, as 
in Example \ref{s6g}. Hence only zero, linear, quadratic or exponential growth is displayed by the cluster recurrences 
(\ref{perrec}). Interestingly, these are the only  types of growth found in the other families of cluster maps 
considered in \cite{gp}. We do not know if other types of growth are possible; are there cluster maps 
with cubic degree growth, for instance?

\section{Poisson and symplectic structures} 

So far we have alluded to the concept of integrability, but have skirted around the issue of giving a precise definition 
of what it means for a map to be integrable. An expected feature of integrability is the ability to find explicit 
solutions of the equations being considered; the recurrence (\ref{a13rec}) displays this feature, because all of 
its iterates satisfy a linear recurrence of the form (\ref{linrec}), which can be solved exactly.   There are many 
other criteria that can be imposed: existence of sufficiently many conserved quantities or symmetries, or compatibility 
of an associated linear system (Lax pair), for instance; and not all of these requirements may be appropriate in 
different circumstances. 
It is an unfortunate fact that the definition of an integrable system 
varies depending on the context, i.e.\ whether it be autonomous or non-autonomous ordinary differential equations, partial differential equations, difference equations, maps 
or something else that is being considered. Thus we need to address  this problem and clarify the context, in order to specify what integrability means for maps associated with cluster algebras. 

There is a precise definition of Liouville integrability in the context of finite-dimensional Hamiltonian mechanics, on a real symplectic 
manifold $M$ of dimension $2m$, with associated Poisson bracket $\{ \,,\,\}$: given a particular function $H$,  the Hamiltonian flow generated 
by $H$ is completely integrable,  in the sense 
of Liouville, if there exist  $m$ independent functions on $M$ (including the Hamiltonian), 
say $H_1=H$, $H_2, \ldots, H_m$, which are in involution with respect to the Poisson bracket, i.e.\ 
$\{H_j,H_k\}=0$ for all $j,k$. In the context of classical mechanics, this notion of integrability 
provides everything one could hope for. To begin with,  systems satisfying 
these requirements have (at least) 
$m$ independent conserved quantities: all of the first integrals $H_1,\ldots, H_m$ are 
preserved by the time evolution, so each of the trajectories lies on an $m$-dimensional intersection of 
 level sets for these functions.  Furthermore,  
Liouville proved that  
the solution of the equations of motion  for such systems 
can be reduced to a finite number of 
quadratures, so they really are ``able to be integrated'' as one would expect; and Arnold 
showed in addition that the flow reduces to quasiperiodic motion on compact $m$-dimensional 
level sets, which are diffeomorphic to tori $T^m$ \cite{arnold}, so nowadays the combined 
result is referred to as the Liouville-Arnold theorem.  Another approach to 
integrability is to require a sufficient number of symmetries, and this is a  consequence of the  Liouville 
definition: the Hamilton's equations arising from $H$ have the maximum number of commuting symmetries, namely the flows 
generated by each of the first integrals $H_j$.   

The notion of Liouville integrability can be extended to symplectic maps in a natural way \cite{bruschi,maeda,veselov}. However, the requirement of working in even dimensions is too restrictive for our purposes, so instead of a symplectic form 
we start with a (possibly degenerate) Poisson structure and consider Poisson maps $\varphi$, 
defined in terms of the pullback of functions, given by $\varphi^* F=F\cdot \varphi$. 

\begin{definition} Given a Poisson bracket $\{\, , \}$ on a manifold $M$, a map $\varphi:\, M\rightarrow M$ is 
called a Poisson map if 
$$ 
\{ \varphi^* F, \varphi^* G\} = \varphi^*\{F,G\}
$$ 
holds for all functions $F,G$ on $M$. 
\end{definition}  

(We are being deliberately vague about what sort of Poisson manifold $(M,\{,\})$  is being considered, e.g.\ a real smooth manifold, or 
a complex algebraic variety, and what sort of functions, e.g.\ smooth/analytic/rational, because this may vary according 
to the context.)  

In order to have a suitable notion of integrability for cluster maps, we first require a 
compatible Poisson structure of some kind. In general, 
given a difference equation or map, there is no canonical way to find a compatible Poisson bracket. Fortunately, it turns 
out that for cluster algebras there is often a natural 
Poisson bracket, of log-canonical type, that is 
compatible with cluster mutations; and there is always a log-canonical presymplectic form \cite{fg,gsv1,gsv2,in}.

\begin{example}\label{s5} 
 {\bf Somos-5 Poisson bracket:}  
The skew-symmetric exchange matrix 
\beq\label{s5b} 
\left(
\begin{array}{ccccc} 
0 & 1 & -1 & -1 & 1 \\
-1 & 0 & 2 & 0 & -1 \\   
1 & -2 & 0 & 2 & -1 \\ 
1 & 0 & -2 & 0 & 1 \\
-1 & 1 & 1 & -1 & 0 
\end{array}\right) 
\eeq 
is cluster mutation-periodic with period 1. Its associated  cluster map is a Somos-5 recurrence, which belongs 
to family (iv) above, given by (\ref{iv}) with $N=5$, $p=1$, $q=2$.   
The skew-symmetric matrix $P=(p_{ij})$ given by 
$$ 
P =  
\left(
\begin{array}{ccccc} 
0 & 1 & 2 & 3 & 4 \\
-1 & 0 & 1 & 2 & 3 \\ 
-2 & -1 & 0 & 1 & 2 \\ 
-3 & -2 & -1 & 0 & 1 \\
-4 & -3 & -2 & -1 & 0 
\end{array}\right) 
$$ 
defines a Poisson bracket, given in terms of the original cluster variables 
${\bf x} = (x_1,x_2,x_3,x_4,x_5)$ by 
\beq\label{s5br} 
\{x_i ,x_j\} = p_{ij} x_i x_j. 
\eeq 
This bracket is called log-canonical because it is just given by the constant matrix $P$  in terms of the logarithmic 
coordinates $\log x_i$. It is also compatible with the cluster algebra structure, in the sense that it remains 
log-canonical under the action of 
any mutation, i.e. writing ${\bf x} \mapsto \mu_k({\bf x}) ={\bf x}'=(x_i')$, in the new cluster variables  it takes 
the form 
$$ 
\{x_i' ,x_j'\} = p'_{ij} x_i' x_j'
$$ 
for some constant skew-symmetric matrix $P'=(p_{ij}')$. Moreover, under the cluster map 
$\varphi=\rho^{-1}\cdot\mu_1$ defined by 
\beq\label{s5map} 
\varphi: \, (x_1,x_2,x_3,x_4,x_5) \mapsto \left(x_2,x_3,x_4,x_5,\frac{x_2x_5+x_3x_4}{x_1}\right), 
\eeq 
the bracket (\ref{s5br}) is preserved, in the sense that for all $i,j\in[1,5]$ the pullback  of the coordinate 
functions 
by the map  
satisfies 
$$ 
\{\varphi^*x_i,\varphi^*x_j\} = \varphi^*\{x_i,x_j\}. 
$$ 
Hence $\varphi$ is a Poisson map with respect to this bracket. 
\end{example}

Given a Poisson map, we can give a definition of discrete integrability, by adapting 
a definition from \cite{pol}, that applies in the continuous case of Hamiltonian  flows on Poisson 
manifolds.

\begin{definition} \label{cint} 
Suppose
that the Poisson tensor is of constant rank $2m$ on a dense open subset of
a Poisson manifold $M$ of dimension $N$, and that the algebra of Casimir
functions is maximal, i.e.\  it contains $N-2m$ independent functions. 
A Poisson map $\varphi: \, M\to M$  is said to be completely integrable if it preserves   
$N-m$ independent functions $F_1,\ldots, F_{N-m}$ which are in
involution, including the 
Casimirs. 
\end{definition} 

\begin{example} {\bf Complete integrability of the $\tilde{A}_{1,3}$ cluster map:} 
Setting $P=B$ with the exchange matrix (\ref{a13b}) in Example \ref{a13}, the bracket 
$$ 
\{x_i,x_j\} =b_{ij}x_ix_j 
$$ 
is compatible with the cluster algebra structure, and is preserved by the cluster map 
$\varphi$ corresponding to (\ref{a13rec}). At points where all coordinates $x_j$ are non-zero, 
the Poisson tensor has full rank 4, since $B$ is invertible; so there are no Casimirs. 
Note that $B^{-1}=-\frac{1}{2}B$, so $B$ is proportional to its own inverse, 
and the map $\varphi$ is symplectic, i.e. $\varphi^*\om = \om$, where up to overall 
rescaling the symplectic form is 
\beq \label{pres}
\om = \sum_{i<j} \frac{b_{ij}}{x_ix_j} \, \rd x_i\wedge \rd x_j. 
\eeq 
Now, observe that the  recurrence can be rewritten  with 
 a $2\times 2$ determinant, as 
$$ 
\left|D_n \right|=1, \quad \mathrm{where} \quad  
D_n = \left( \begin{array}{cc} 
x_n & x_{n+1} \\ 
x_{n+3} & x_{n+4} \end{array} \right), 
$$ 
and construct the $3\times 3$  matrix sequence 
$$ 
\tilde{D}_n =\left(\begin{array}{ccc} 
x_n & x_{n+1} & x_{n+2} \\ 
x_{n+3} & x_{n+4} & x_{n+5} \\ 
x_{n+6} & x_{n+7} & x_{n+8} \end{array} 
\right).  
$$ 
Then, by the method of Dodgson condensation \cite{dodgson}, the determinant can be 
expanded as 
$$ 
| \tilde{D}_n | =\frac{1}{x_{n+4}} \left| \begin{array}{cc} \left| D_n\right| &  \left| D_{n+1}\right| \\ 
 \left| D_{n+3}\right| &  \left| D_{n+4} \right| \end{array} \right| = 0 . 
$$     
Further calculation shows that the kernel of $\tilde{D}_n$ is spanned 
by the vector $(1,-J_n,1)^T$, where $J_n$  is periodic with period 3, so there is a linear relation 
\beq\label{jrec} 
x_{n+2}-J_n x_{n+1}+x_n =0, \quad \mathrm{with} \quad J_{n+3}=J_n. 
\eeq 
Similarly, the linear relation (\ref{linrec}), with invariant $K$ (independent of $n$), corresponds to 
the fact that the kernel of $\tilde{D}_n^T$ is spanned by $(1,-K,1)^T$.   The $J_i$ 
can be considered as functions of the phase space coordinates $x_1$, $x_2$, $x_3$, $x_4$, by writing 
$$ 
J_1=\frac{x_1+x_3}{x_2}, 
$$ 
and similarly for $J_2,J_3$. Computing the Poisson bracket between these functions yields 
$$ 
\{J_i,J_{i+1}\} = 2J_iJ_{i+1}-2, 
$$ 
where the indices are read $\bmod\,3$, so that $J_1,J_2,J_3$ form a Poisson subalgebra of 
dimension 3; and  
$$ 
K=J_1J_2J_3 -J_1-J_2-J_3
$$ 
is a Casimir for this subalgebra, in the sense that $\{J_i,K\}=0$ for $i=1,2,3$. 
The map $\varphi$ preserves any symmetric function of the $J_i$, so picking the 
three independent functions 
$$ 
F_1=K, \quad F_2=J_1J_2+J_2J_3+J_3J_1, \quad F_3= J_1+J_2+J_3 
$$ 
we have 
$\varphi^*F_j=F_j$ for all $j$, but at most two of these can be in  involution: 
$\{F_1,F_2\}=0=\{F_1,F_3\}$, but $\{F_2,F_3\}\neq 0$. Thus, choosing just $F_1$ and $F_2$, say, 
the conditions of Definition \ref{cint} are satisfied, and the map $\varphi$ given by (\ref{a13map}) 
 is completely integrable.  
\end{example} 

\begin{remark} The fact that  cluster variables  obtained from affine quivers satisfy linear relations with constant 
coefficients, such as (\ref{linrec}),  has been shown 
in various different ways: for type $A$  in \cite{fm,FH}, using Dodgson condensation (equivalently, 
the Desnanot-Jacobi formula);  for types $A$ and $D$ in the context  of frieze relations \cite{ars}; and for all simply-laced 
types $A$, $D$, $E$ in \cite{ks}, using cluster categories (but see also \cite{difk,pasha} for another family of quivers 
made from  products of finite and affine Dynkin types $A$). The fact that there are additional  linear relations with periodic coefficients, like (\ref{jrec}), 
was shown for all $\tilde{A}_{p,q}$ quivers in \cite{FH}, where 
it was also found that the quantities $J_i$ are coordinates in the dressing chain for Schr\"odinger 
operators, and this has recently been extended to affine types $D$ and $E$ \cite{joe}. 
\end{remark} 

\begin{example} \label{a12}  
{\bf Non-existence of a log-canonical bracket for $\tilde{A}_{1,2}$:} 
For the cluster algebra of type $\tilde{A}_{1,2}$, defined by the skew-symmetric exchange matrix 
$$ 
B=  \left(\begin{array}{ccc} 
0 & 1 & 1 \\ 
-1 & 0 & 1 \\ 
-1 & -1 & 0 \end{array}\right) 
$$ 
it is easy to verify that there is no bracket of log-canonical form, like (\ref{s5br}), that is compatible 
with cluster mutations. However, iterates of the cluster map, defined by the recurrence  
$$ 
x_nx_{n+3}=x_{n+1}x_{n+2}+1, 
$$ 
satisfy the linear relation 
$
x_{n+4}-Kx_{n+2}+x_n=0$, 
for a first integral $K$. In fact, setting $u_n=x_nx_{n+1}$ yields a recurrence of second order, 
\beq\label{a12u} 
u_nu_{n+2} =u_{n+1}(u_{n+1}+1), 
\eeq 
and, rewriting $K$ in terms of $u_1,u_2$, this corresponds to a symplectic map $\hat\varphi$ in the 
$(u_1,u_2)$ plane with symplectic form $\hat\om=\rd \log u_1\wedge  \rd \log u_2$ and one first integral; so the 
map $\hat\varphi$ is completely integrable.   
\end{example} 

\begin{example} {\bf Casimirs for Somos-5:} The Poisson tensor for the Somos-5 map (\ref{s5map}), defined 
by (\ref{s5br}), has rank 2 on $\C^5\setminus\{x_i=0\}$ (away from the coordinate hyperplanes). 
The kernel of the matrix $P$ is spanned 
by the vectors 
\beq\label{pal} 
\tilde{{\bf v}}_1=(1,-2,1,0,0)^T, \, 
\tilde{{\bf v}}_2=(0,1,-2,1,0)^T, \, 
\tilde{{\bf v}}_3=(0,0,1,-2,1)^T, 
\eeq 
which  correspond to three independent Casimir functions
$$ 
F_1=\bx^{\tilde{{\bf v}}_1}=\frac{x_1x_3}{x_2^2}, \quad
F_2=\bx^{\tilde{{\bf v}}_2}=\frac{x_2x_4}{x_3^2}, \quad
F_3=\bx^{\tilde{{\bf v}}_3}=\frac{x_3x_5}{x_4^2}, 
$$ 
whose Poisson bracket with any other function $G$ vanishes: $\{F_j,G\}=0$ for $j=1,2,3$. 
There are  two independent first integrals $H_1$, $H_2$, i.e.\ functions that are preserved 
by the action of $\varphi$, so that  $\varphi^*H_i = H_i\cdot \varphi=H_i$ for $i=1,2$; and these 
are themselves Casimirs because they can be written in terms of the $F_j$ \cite{honetams}: 
\beq\label{h1s5}
H_1 = F_1F_2F_3+\frac{1}{F_1}+ \frac{1}{F_2}+\frac{1}{F_3} +\frac{1}{F_1F_2F_3} ,
\eeq 
\beq\label{h2s5} 
H_2 = F_1F_2+F_2F_3+\frac{1}{F_1F_2}+\frac{1}{F_2F_3} +\frac{1}{F_1F_2^2F_3} .
\eeq 
However, the full algebra of Casimirs is not preserved by the map $\varphi$, because $F_1,F_2,F_3$  transform as 
$$
\varphi^*F_1=F_2, \quad \varphi^*F_2=F_3, \quad \varphi^*F_3=\frac{F_2F_3+1}{F_1F_2^2F_3^2}.
$$
 Hence the Somos-5 map is not completely integrable with respect to this bracket. 
\end{example} 

The previous two examples show that if the exchange matrix $B$ is degenerate, then  the cluster coordinates may not be the correct ones to use, as either there 
is no invariant log-canonical  bracket in these coordinates, as in the case of $\tilde{A}_{1,2}$,  or even if there is 
such a bracket, a 
full set of Casimirs is not preserved by the cluster map. (A Poisson map sends Casimirs to other Casimirs, but need 
not preserve each Casimir individually.) The way out of this quandary, which was already hinted at in Example \ref{a12}, 
is to work on a reduced space where the map $\varphi$ reduces to a symplectic map $\hat\varphi$. It turns out that there 
is a canonical way to do this, based on the presymplectic form $\om$ associated with the cluster algebra, which in general, 
for any skew-symmetric exchange matrix $B=(b_{ij})$, is given  by the formula (\ref{pres}) above. 

In the case that $B$ is nondegenerate (which is possible for even $N$ only, as in Example \ref{a13}), $\om$ is a closed, nondegenerate 2-form, so 
the cluster map is symplectic, but otherwise $\om$ has a null distribution, generated  by vector fields of the form 
$$ 
\sum_{j=1}^N w_j x_j \frac{\partial}{\partial x_j}, \quad \mathrm{for} \quad {\bf w}=(w_j)\in\mathrm{ker} \, B.
$$  
These vector fields all commute with other, and can be integrated to yield  a commuting set of scaling symmetries: 
each ${\bf w}\in\mathrm{ker}\, B$ generates a one-parameter scaling group 
\beq\label{sym} 
{\bf x}\mapsto \tilde{\bx} = \la^{\bf w}\cdot \bx, \qquad \la \in \C^*, 
\eeq  
where the notation means that each component is scaled so that $\tilde{x}_j=\la^{w_j}x_j$. 
Regarding $B$ as a linear transformation on $\Q^N$, skew-symmetry means that there is an orthogonal direct 
sum decomposition $\Q^N = \mathrm{im}\, B \oplus \mathrm{ker} \, B$. 
If $B$ has rank $2m$, then 
an integer  basis ${\bf v}_1,{\bf v}_2,\ldots, {\bf v}_{2m}$ for 
$\mathrm{im}\, B$ yields a complete set of rational functions invariant under the symmetries (\ref{sym}),  given 
by the monomials 
\beq\label{mon}
u_j = \bx^{{\bf v}_j}, \quad j=1,\ldots, 2m.
\eeq 
In the case that $B$ has period 1, it was shown in \cite{FH} that  by choosing the basis suitably, the rational  map 
$
\pi: \, \bx \mapsto {\bf u}=(u_j) 
$
reduces $\varphi$ to a birational symplectic map $\hat\varphi$ in dimension $2m$, with symplectic form 
$\hat\om$, in the sense that 
$\hat{\varphi}\cdot \pi=\pi\cdot \varphi$, and $\pi^*\hat\om = \om$, 
where 
\beq\label{symp}
\hat\om = \sum_{i<j} \frac{\hat{b}_{ij}}{u_iu_j} \, \rd u_i \wedge \rd u_j 
\eeq 
(for a certain skew-symmetric matrix $\hat{B}=(\hat{b}_{ij})$)
is also log-canonical. In \cite{hi} it was further shown that (up to an overall sign) there is a canonical choice of basis for 
$\mathrm{im}\, B   \cap \Z^N$ with the property that 
$$ 
\varphi^* \bx^{{\bf v}_j}=\bx^{{\bf v}_{j+1}}, \quad \mathrm{for} \quad j\in [1,2m-1].
$$
This is called a {\it palindromic basis}, because the first $N-2m+1$ entries of ${\bf v}_1$ form a 
palindrome, with the remaining $2m-1$ entries being zero, and this palindrome is just shifted along to get the 
other basis elements; the basis 
is fixed uniquely if the first entry of ${\bf v}_1$ is chosen to be positive. 
The advantage of 
a palindromic basis is that the birational  map $\hat\varphi$ is equivalent to an iteration of 
a single recurrence relation. 

\begin{definition} Given a cluster mutation-periodic skew-symmetric 
exchange matrix $B$ with period 1, of rank $2m$,  and the symplectic coordinates 
$(u_j)\in \C^{2m}$ defined by (\ref{mon}) with  a palindromic basis, the {\it U-system} is the  recurrence corresponding to the reduced cluster 
map $\hat\varphi$, which, for some rational function ${\cal F}$, has  the form 
\beq\label{usys}
u_n u_{n+2m} = {\cal F}(u_{n+1},\ldots, u_{n+2m-1}).  
\eeq 
\end{definition}   

We have already seen an example of a U-system, namely the reduced  recurrence (\ref{a12u}) for $\tilde{A}_{1,2}$. 
An integrable U-system corresponds to the canonical version of integrability for maps: 
the U-system   is equivalent to a symplectic map in dimension $2m$, so $m$ independent first integrals in involution 
are needed for complete integrability. 

\begin{example} {\bf Complete integrability of the Somos-5 U-system:} With $B$ given by (\ref{s5b}), 
a palindromic basis for $\mathrm{im}\, B$ is written using 
 (\ref{pal}) as  
$$ 
{\bf v}_1 = \tilde{{\bf v}}_1+   \tilde{{\bf v}}_2 = (1,-1,-1,1,0)^T, \quad 
{\bf v}_2 = \tilde{{\bf v}}_2+   \tilde{{\bf v}}_3 = (0, 1,-1,-1,1)^T, 
$$ 
so the reduced coordinates are 
$$ 
u_1=\frac{x_1x_4}{x_2x_3}=F_1F_2, \quad 
u_2=\frac{x_2x_5}{x_3x_4}=F_2F_3, 
$$ 
and 
$\om = \sum_{i<j}b_{ij}\rd\log x_i\wedge \rd \log x_j$ reduces to the 
symplectic form 
$$ 
\hat\om=\frac{\rd u_1\wedge \rd u_2}{u_1u_2} 
$$ in these coordinates. 
The cluster map (\ref{s5map}) reduces to an iteration  of the U-system 
\beq\label{s5u}
u_nu_{n+2}=\frac{u_{n+1}+1}{u_{n+1}},  
\eeq 
and although $H_1$ does not survive this reduction, the first integral $H_2$ can be rewritten 
in terms of $u_1,u_2$, to yield the 
function   
$$ 
H =u_1 +u_2 + \frac{1}{u_1}+ \frac{1}{u_2} + \frac{1}{u_1u_2}, 
$$ 
so the U-system corresponds to a completely integrable symplectic map in two 
dimensions. The generic level sets of $H$ are cubic curves of genus 1, and this is an 
example of a symmetric QRT map (see \cite{hhkq} and references).  
\end{example} 

The general Somos-6 recurrence, with constant coefficients $\al, \bet, \gamma$,  has the form 
\beq\label{s6gen} 
x_n x_{n+6} = \al x_{n+1} x_{n+5} +\bet x_{n+2}x_{n+4} + \gam x_{n+3}^2, 
\eeq 
which has the Laurent property \cite{fz}, but cannot come from a cluster algebra when 
$\al\bet\gamma\neq 0$, due to there being too many terms on the right-hand side. 
In fact, it appears in the more general setting of mutations in LP algebras, which allow 
exchange relations with more terms \cite{lp}.  Being quadratic relations, Somos recurrences 
are reminiscent of Hirota bilinear equations for tau functions in soliton theory, and indeed, 
the general Somos-6 recurrence is a reduction of Miwa's equation \cite{Date}, which is the bilinear 
discrete BKP equation, also known as the cube recurrence in algebraic combinatorics. 
Here we conclude our discussion of Example \ref{s6}, by setting $\beta=0$, to obtain 
a bilinear equation with a total of three terms, which can be obtained as a reduction of the 
discrete Hirota equation (bilinear discrete KP, or octahedron recurrence), that is 
\beq \label{dkp} 
T_1T_{-1}=T_2T_{-2}+T_3T_{-3}, 
\eeq   
where the tau function $T=T(m_1,m_2,m_3)$ and the subscript $\pm j$ denotes a shift 
in the $j$th independent variable, so e.g.\  $T_{\pm 1} = T(m_1\pm 1,m_2,m_3)$, and 
so on. The advantage of making a reduction from this equation with more independent 
variables is that  it has a Lax pair, which  reduces to a Lax pair 
for the Somos recurrence, and there is an associated spectral curve, whose 
coefficients provide 
first integrals.  

\begin{example}{\bf A Somos-6 U-system:}
Setting $\bet=0$
in (\ref{s6gen}) produces 
\beq\label{s6n} 
x_n x_{n+6} = \al x_{n+1} x_{n+5} + \gam x_{n+3}^2.
\eeq 
This differs from (\ref{s6rec}) and (\ref{s6map}) by the inclusion of coefficients $\al,\gamma$, 
which can be achieved by augmenting the cluster algebra with frozen variables that 
appear in the exchange relations but do not themselves mutate (see \cite{fm} and references, 
for instance), and does not change other features such as Poisson brackets or the (pre)symplectic forms. 
Upon applying the method in \cite{hkw}, we can obtain (\ref{s6gen}) as a plane wave reduction 
of (\ref{dkp}), by setting 
$$ 
T(m_1,m_2,m_3)=a_1^{m_1^2}a_2^{m_2^2}a_3^{m_3^2}\, x_n, \qquad n=m_0 +   3m_1 + 2m_3, 
$$ 
with $m_0$ arbitary, and taking $\al =a^2_3/a^2_1$, $\gamma=a^2_2/a^2_1$. 
Under this reduction, the linear system whose compatibility gives the discrete KP equation becomes 
\beq\label{linsys} 
\begin{array}{rcl} 
Y_n \psi_{n+3} + \alpha \ze \psi_{n+2} & = & \xi \psi_n, \\ 
\psi_{n+3} -X_n \psi_{n+1}  & = & \ze \psi_n, 
\end{array} 
\eeq 
where $\psi_n$ is a wave function, $\ze,\xi$ are spectral parameters, 
and  
$$ 
X_n = \frac{x_{n+2}x_{n+3}}{x_{n+4}x_{n+1}}, \quad 
Y_n = \frac{x_{n+4}x_n}{x_{n+3}x_{n+1}}. 
$$ 
The equation (\ref{s6n}) is the compatibility condition for these two linear equations for $\psi_n$
(to be precise, the parameter $\gamma$ arises as an integration constant). This is more conveniently seen 
by writing the second linear equation in matrix form, with a vector $\Psi_n=(\psi_n,\psi_{n+1},\psi_{n+2})^T$, 
as 
\beq\label{meqn}
\Psi_{n+1}={\bf M}_n \psi_n, \quad {\bf M}_n =\left(\begin{array}{ccc} 
0 & 1 & 0 \\ 
0 & 0 & 1 \\ 
\ze & X_n & 0 \end{array}\right),
\eeq  
and then using the second linear  equation in (\ref{linsys})  to reformulate the first one as an 
eigenvalue problem with $\Psi_n$ as the eigenvector, that is 
\beq\label{leqn}
{\bf L}_n \Psi_n=\xi\Psi_n, \quad {\bf L}_n =\left(\begin{array}{ccc} 
\ze Y_n & u_n & \ze\al \\ 
\ze^2\al & \ze(Y_{n+1}+\al X_n) & u_{n+1} \\ 
\ze u_{n+2} & \ze^2 \al + u_{n+1}^{-1} & \ze (Y_{n+2}+\al X_{n+1}) \end{array}\right).
\eeq  
In the above expression for the Lax matrix ${\bf L}_n$, we have introduced the quantities 
$$ 
u_n = X_nY_n = \frac{x_nx_{n+2}}{x_{n+1}^2}, 
$$ 
which for $n=1,2,3,4$ give a set of symplectic coordinates obtained from the palindromic 
basis 
${\bf v}_1=(1,-2,1,0,0,0)^T$,  ${\bf v}_2=(0,1,-2,1,0,0)^T$,  
${\bf v}_3=(0,0,1,-2,1,0)^T$,
${\bf v}_4=(0,0,0,1,-2,1)^T$ for $\mathrm{im}\, B$, with $B$ as in (\ref{s6b}), and satisfy  
the U-system 
\beq \label{us6}
u_nu_{n+4}=\frac{\al u_{n+1}u_{n+2}^2u_{n+3}+\gamma} {u_{n+1}^2u_{n+2}^3u_{n+3}^2}
\eeq   
(which should be compared with the tropical formulae (\ref{utrop}) and (\ref{usystrop}) above), 
correponding to the reduced cluster map $\hat\varphi$.  
The symplectic form $\hat \om$, such that $\hat\varphi^*\hat\om=\hat\om$,  is 
$$ 
\hat\om=\sum_{i<j}\hat{b}_{ij} \rd \log u_i \wedge \rd\log u_j, \quad 
\hat{B}=(\hat{b}_{ij}) = \left(\begin{array}{cccc} 0 & 1 & 2 & 1 \\
-1 & 0 & 2 & 2 \\
-2 & -2 & 0 & 1 \\ 
-1 & -2 & -1 & 0 \end{array}\right), 
$$ 
so the associated nondegenerate Poisson bracket for these coordinates 
is given by  $\{u_i,u_j\}=\hat{p}_{ij}u_iu_j$ with $(  \hat{p}_{ij})=\hat{B}^{-1}$. 
The compatibility condition of the matrix system given by (\ref{meqn}) and (\ref{leqn}) 
is the discrete Lax equation 
$$ 
{\bf L}_{n+1} {\bf M}_n ={\bf M}_n {\bf L}_n ,
$$ 
which is equivalent to the U-system (\ref{us6}). So this is an isospectral evolution, 
and the spectral curve 
\beq \label{spec}
\det ({\bf L}_n -\xi \mathbf{1}) = -\xi^3+H_1\xi^2+(1-H_2\ze^2)\xi+\al^3\ze^5+\gamma\ze^3=0 
\eeq
is independent of $n$, with the non-trivial coefficients being $H_1$, given by 
$$ 
(u_n+u_{n+3})u_{n+1}u_{n+2} 
+\al\left(\frac{1}{u_{n}u_{n+1}}+ \frac{1}{u_{n+1}u_{n+2}}+\frac{1}{u_{n+2}u_{n+3}}\right)
+ \frac{\gamma}{u_nu_{n+1}^2u_{n+2}^2u_{n+3}}, 
$$ 
and 
$$ 
\begin{array}{rcl} 
H_2 & = & \al\left(\frac{u_nu_{n+1}}{u_{n+3}}+\frac{u_{n+2}u_{n+3}}{u_{n}} \right)
+\gamma\left(\frac{1}{u_nu_{n+1}u_{n+2}}+ \frac{1}{u_{n+1}u_{n+2}u_{n+3}}\right) \\
&&  
+ \al^2 \left(\frac{1}{u_nu_{n+1}^2u_{n+2}}+ \frac{1}{u_{n+1}u_{n+2}^2u_{n+3}}\right) 
+\frac{\al\gamma }{u_nu_{n+1}^3u_{n+2}^3u_{n+3}}, 
\end{array}
$$ 
which provide  two independent first integrals. It can be verified directly that $\{H_1,H_2\}=0$, 
which shows that each iteration of (\ref{us6}) corresponds to a completely integrable 
symplectic map $\hat\varphi$ (in different coordinates, the involutivity of these quantities was also shown in \cite{hones6}). 
The trigonal spectral curve (\ref{spec}) has genus 4, and admits the involution 
$(\ze,\xi)\mapsto(-\ze,-\xi)$, giving a quotient curve of genus 2, with a Prym variety 
that is isomorphic to the Jacobian of a second genus 2 curve, analogous to the situation 
for the general Somos-6 map in \cite{fh}. However, in this case there  is a more direct way to find the second 
genus 2 curve, as the hyperelliptic spectral curve of a $2\times 2$ Lax pair 
obtained by deriving (\ref{s6n}) as a reduction 
of a discrete time Toda  equation on a 5-point lattice \cite{hkw,hkq}. 
For explicit analytic formulae for the solutions in terms of genus 2 sigma functions, see 
\cite{hones6,fh}.  
\end {example}

\section{Discrete Painlev\'e equations from  coefficient mutation}

The continuous Painlev\'e equations are a special set of non-autonomous 
ordinary differential equations 
of second order  that are characterized by the absence of movable critical points in their 
solutions, which is known as the Painlev\'e property. 
Discrete Painlev\'e equations are a particular class of ordinary difference equations 
which, like their continuous counterparts, 
are non-autonomous 
(meaning that the independent variable appears explicitly); in many cases, they appeared  
from the search for an appropriate discrete 
analogue of the Painlev\'e property \cite{grp}. The resulting notion of 
singularity confinement turned out to be much weaker than the Painlev\'e 
property for differential equations, and is not sufficient 
for integrability, although it is a very useful tool when used judiciously 
in tandem with other techniques for identifying integrable maps or discrete Painlev\'e equations \cite{con2}.
In fact, singularity confinement seems to be very closely related to the Laurent property \cite{honeconf}, 
and it is interesting to speculate whether all discrete integrable systems are related to a system with 
the Laurent property by  introducing a tau function or some other lift of the coordinates \cite{chang,hk,mase1}. 

Recently there have been various studies that show how certain discrete Painlev\'e equations and their higher 
order analogues can arise from mutation of coefficients in cluster algebras. Here we concentrate on the methods 
used in \cite{hi}, but for other related approaches see the work of Okubo \cite{okubo,os} and that of Bershtein et al.\ \cite{bgm}. 

A {\it Y-system} is  a set of difference equations arising  as relations between coefficients appearing from a sequence of mutations
in a cluster algebra with periodicity. The original Y-systems were obtained by Zamolodchikov as a set of functional equations 
in certain quantum field theories associated with  simply-laced affine Lie algebras \cite{zam}, yet they arise from cluster algebras 
of finite type obtained from the corresponding finite-dimensional root systems, and display purely periodic dynamics. 
Generalized Y-systems were defined by Nakanishi \cite{nak} starting from a general notion of periodicity in a cluster algebra, 
and typically display  complicated dynamical behaviour. 

Here we concentrate on the case 
of cluster mutation-periodic quivers with period 1, for which the Y-system can be written as a single scalar 
difference equation, given by 
\beq\label{ysys} 
y_ny_{n+N} = \frac{ \prod_{j=1}^{N-1} (1+y_{n+j})^{[a_j]_+}   } 
{\prod_{j=1}^{N-1} (1+y_{n+j}^{-1})^{[-a_j]_+}},
\eeq 
where, as in Theorem \ref{peri},    $a_j=b_{1,j+1}$ are the components of the palindromic 
$(N-1)$-tuple that determines the exchange matrix. (Here 
we assume that the first non-zero component $a_j$ is positive; there is no loss of generality in doing so, 
 due to the freedom to replace $B\to -B$, but some signs are reversed compared with \cite{hi} and \cite{nak}.)
In this context, the coefficient-free recurrence 
(\ref{perrec}) that defines 
the cluster map  is referred to as the {\it T-system}. It was first  observed in \cite{f1} that there is a relation between 
the evolution of coefficients ${\bf y}$ under mutations (\ref{comut})  in a cluster algebra,  and the evolution of 
cluster variables ${\bf x}$ due to the    associated coefficient-free cluster mutations given by (\ref{freeclmut}), which can 
be summarized by the slogan that ``the T-system provides a solution of the Y-system.'' In the 
case at hand, the precise statement is that making  the subsitution  
\beq\label{ysub} 
y_n = \prod_{j=1}^{N-1}x_{n+j}^{a_j}  
\eeq    
in (\ref{ysys}) provides a solution of the Y-system whenever $x_n$ satisfies the coefficient-free T-system 
(\ref{perrec}). 
 
Although the equations (\ref{perrec}) and (\ref{ysys}) are both 
of order $N$, there can be a discrepancy between the solutions of the T-system and the Y-system, in the sense 
that the general solution of the former does not yield the general solution of the latter. This discrepancy is 
determined by the following result. 

\begin{proposition} Let $x_n$ satisfy the modified T-system 
\beq\label{tz} 
x_nx_{n+N}={\cal Z}_n \left(\prod_{j=1}^{N-1}x_{n+j}^{[a_j]_+} + \prod_{j=1}^{N-1}x_{n+j}^{[-a_j]_+}
\right). 
\eeq 
Then the substitution (\ref{ysub}) yields a solution of the Y-system (\ref{ysys}) if and only if ${\cal Z}_n$ 
satisfies the {\it Z-system} 
\beq\label{z} 
\prod_{j=1}^{N-1} {\cal Z}_{n+j}^{a_j}=1. 
\eeq 
\end{proposition} 

Each iteration of the modified T-system (\ref{tz}) with non-autonomous coefficients evolving according to (\ref{z}) 
preserves the presymplectic form given by (\ref{pres}) in terms of the entries of the exchange matrix $B$, and if $B$ is 
degenerate we can use  a 
palindromic basis for $\mathrm{im}\, B$  to reduce this to a non-autonomous recurrence in lower dimension that 
preserves the symplectic form (\ref{symp}). 
 
\begin{definition} The pair of equations (\ref{tz}) and (\ref{z}) is called the {\it $T_z$-system}. The {\it $U_z$-system} 
associated with (\ref{tz}) is given by (\ref{z}) together with 
$$ 
u_n u_{n+2m} = {\cal Z}_n\,{\cal F}(u_{n+1},\ldots, u_{n+2m-1}), 
$$ 
where the rational function $\cal F$ is the same as in (\ref{usys}). 
\end{definition} 

We conclude this section with a couple of examples.

\begin{example} {\bf Somos-5 Y-system and q-Painlev\'e II:} 
The Y-system associated with the exchange matrix (\ref{s5b}) is 
$$ 
y_ny_{n+5}=\frac{(1+y_{n+1})(1+y_{n+4})}{(1+y_{n+2}^{-1})(1+y_{n+3}^{-1})}, 
$$ 
and (noting that on the 
right-hand side of the substitution (\ref{ysub}) there is the freedom to shift $n\to n+1$) 
the general solution of this can be written as 
$$
y_n =\frac{x_{n}x_{n+3}}{x_{n+1}x_{n+2}}, 
$$  
where $x_n$ satisfies the non-autonomous Somos-5 relation  
$$ 
x_nx_{n+5}={\cal Z}_n(x_{n+1}x_{n+4}+x_{n+2}x_{n+3}), \quad \mathrm{with} \quad 
\frac{{\cal Z}_n{\cal Z}_{n+3}}{{\cal Z}_{n+1}{\cal Z}_{n+2}}=1.
$$ 
Equivalently, we can identify $y_n=u_n$ and solve the third order Z-system for ${\cal Z}_n$ to 
write the $U_z$-system as a non-autonomous version of the QRT map (\ref{s5u}), that is 
$$
u_nu_{n+2} = {\cal Z}_n (1+u_n^{-1}), \quad \mathrm{with} \quad {\cal Z}_n=\bet_n {\mathfrak q}^n, \,\, 
\bet_{n+2}=\bet_n.
$$ 
The latter is equivalent to 
a q-Painlev\'e II equation identified in \cite{kruskal}, having a continuum limit to the 
Painlev\'e II differential equation 
$$
\frac{d^2u}{dz^2}=2u^3+ zu+\al.
$$ 
\end{example}
\begin{example}{\bf A q-Somos-6 relation:} The Y-system corresponding to (\ref{s6b}) 
is 
$$ 
y_ny_{n+6}=\frac{(1+y_{n+1})(1+y_{n+5})}{(1+y_{n+3}^{-1})^2}.
$$ 
Its general solution can be written as  
$$ 
y_n =\frac{x_n x_{n+4}}{y_{n+2}^2},  
$$ 
where $x_n$ satisfies a q-Somos-6 relation given by 
\beq\label{qs6}
x_nx_{n+6}= {\cal Z}_n 
(x_{n+1}x_{n+5}+x_{n+3}^2), \quad \mathrm{with} \quad {\cal Z}_n =\alpha_{\pm}{\mathfrak q}_{\pm}^n, 
\eeq 
with the solution of the fourth order Z-system 
$$ 
\frac{{\cal Z}_n {\cal Z}_{n+4}}{{\cal Z}_{n+2}^2}=1 
$$ 
being given in terms of quantities $\al_\pm$ and ${\mathfrak q}_\pm$ that alternate with the parity of $n$. 
Alternatively, one can write 
$$ y_n = u_n u_{n+1}^2u_{n+2}$$ 
with $u_n$  satisfying a non-autonomous version of (\ref{us6}), that is 
$$ 
u_nu_{n+4}=\frac{{\cal Z}_n(u_{n+1}u_{n+2}^2u_{n+3}+1)}{u_{n+1}^2u_{n+2}^3u_{n+3}^2}, 
$$ 
with ${\cal Z}_n$ as in (\ref{qs6}). The latter should be regarded as a fourth order analogue of  a discrete 
Painlev\'e equation.  
\end{example} 

\section{Conclusions} 
We have just scratched the surface in this brief  introduction to cluster algebras and discrete integrability. 
Among  other important examples that we have not described here, 
we would like to mention pentagram maps \cite{gstv} and cluster integrable systems 
related to dimer models \cite{eager,gk}. A slightly different viewpoint, with some different choices of topics, can be found in the review \cite{grup}. 

\section*{Acknowledgements}
The work of all three authors is supported by EPSRC Fellowship EP/M004333/1. ANWH is grateful to the School 
of Mathematics \& Statistics, UNSW for hospitality and additional support under the Distinguished Researcher Visitor Scheme.  
\strut\hfill

\end{document}